\RequirePackage{fix-cm}
\documentclass[smallextended]{svjour3}       

\smartqed  
\usepackage{graphicx}
\usepackage{amssymb}
\usepackage{mathrsfs}
\usepackage{amsmath}
\usepackage{amsfonts}
\usepackage{amssymb}
\usepackage{graphicx} 
\usepackage{xspace} 
\usepackage{bm} 
\usepackage{threeparttable} 
\usepackage{subfigure}
\usepackage{stmaryrd} 
\usepackage[ruled,lined,linesnumbered]{algorithm2e}
\usepackage[normalem]{ulem}
\usepackage{todonotes}
\usepackage{hyperref}
\hypersetup{
  colorlinks=true,
  linkcolor=blue,
  citecolor=blue,
  urlcolor=blue
}
  
\newcommand{\fnc}[1]{\ensuremath{\mathcal{#1}}}
\newcommand{\bfnc}[1]{\ensuremath{\bm{\fnc{#1}}}} 

\newcommand{\mat}[1]{\ensuremath{\mathsf{#1}}}

\newcommand{\M}[0]{\mat{H}}

\newcommand{\Hk}[0]{\mat{H}_k}
\newcommand{\Dxk}[0]{\mat{D}_{x,k}}
\newcommand{\Dyk}[0]{\mat{D}_{y,k}}
\newcommand{\Qxk}[0]{\mat{Q}_{x,k}}
\newcommand{\Qyk}[0]{\mat{Q}_{y,k}}
\newcommand{\Sxk}[0]{\mat{S}_{x,k}}

\newcommand{\Exk}[0]{\mat{E}_{x,k}}
\newcommand{\Eyk}[0]{\mat{E}_{y,k}}

\newcommand{\uk}[0]{\bm{u}_k}
\newcommand{\vk}[0]{\bm{v}_k}

\newcommand{\barMk}[0]{\overbar{\M}_{k}}
\newcommand{\overbar}[1]{\mkern 1.5mu\overline{\mkern-1.5mu#1\mkern-1.5mu}\mkern 1.5mu}

\newcommand{\wk}[0]{\bm{w}_{k}}

\newcommand{\Pk}{\mat{P}_{k}}
\newcommand{\barPk}{\overbar{\mat{P}}_{k}}
\newcommand{\tM}[0]{\tilde{\mat{H}}}
\newcommand{\tDx}[0]{\tilde{\mat{D}}_{x}}
\newcommand{\tQx}[0]{\tilde{\mat{Q}}_{x}}
\newcommand{\tEx}[0]{\tilde{\mat{E}}_{x}}
\newcommand{\tQy}[0]{\tilde{\mat{Q}}_{y}}
\newcommand{\tEy}[0]{\tilde{\mat{E}}_{y}}

\newcommand{\ugd}{\tilde{\bm{u}}}
\newcommand{\vgd}{\tilde{\bm{v}}}
\newcommand{\wgd}{\tilde{\bm{w}}}

\newcommand{\Vtilde}[0]{\tilde{\mat{V}}_k}

\DeclareMathOperator{\spn}{span}

\newcommand{\rank}[0]{\operatorname{rank}}

\newtheorem{thrm}{Theorem}

\newcommand{\etal}[0]{{\em et~al.\@}\xspace}

\newcommand{\ie}[0]{{i.e.\@}\xspace}
\newcommand{\viz}[0]{{viz.\@}\xspace}

\newcommand{\ignore}[1]{} 

\journalname{Journal of Scientific Computing}
\begin{document}

\title{Entropy-stable discontinuous Galerkin difference methods for hyperbolic conservation laws
}


\author{Ge Yan         \and
        Sharanjeet Kaur \and
        Jeffery W. Banks \and
        Jason E. Hicken
}


\institute{Ge Yan, Graduate Student \at
              Rensselaer Polytechnic Institute\\
              \email{geyan3566@gmail.com}           
          \and
              Sharanjeet Kaur, Graduate Student \at
              Rensselaer Polytechnic Institute\\
              \email{kaurs3@rpi.edu}
          \and
              Jeffrey W. Banks, Associate Professor \at
              Rensselaer Polytechnic Institute \\
              \email{banksj3@rpi.edu}
          \and
            Jason E. Hicken, Associate Professor \at
            Rensselaer Polytechnic Institute\\
            \email{hickej2@rpi.edu}
}

\date{Received: date / Accepted: date}

\maketitle

\begin{abstract}
The paper describes the construction of entropy-stable discontinuous Galerkin difference (DGD) discretizations for hyperbolic conservation laws on unstructured grids.  The construction takes advantage of existing theory for entropy-stable summation-by-parts (SBP) discretizations.  In particular, the paper shows how DGD discretizations --- both linear and nonlinear --- can be constructed by defining the  SBP trial and test functions in terms of interpolated DGD degrees of freedom.  In the case of entropy-stable discretizations, the entropy variables rather than the conservative variables must be interpolated to the SBP nodes.  A fully-discrete entropy-stable scheme is obtained by adopting the relaxation Runge-Kutta version of the midpoint method.  In addition, DGD matrix operators for the first derivative are shown to be dense-norm SBP operators.  Numerical results are presented to verify the accuracy and entropy-stability of the DGD discretization in the context of the Euler equations.  The results suggest that  DGD and SBP solution errors are similar for the same number of degrees of freedom.  Finally, an investigation of the DGD spectra shows that spectral radius is relatively insensitive to discretization order; however, the high-order methods do suffer from the linear instability reported for other entropy-stable discretizations.
\end{abstract}
\section*{Mathematics Subject Classification}
65M60, 65M70, 65M12

\section{Introduction}
\label{intro}
Many studies have demonstrated that high-order discretizations can simulate hyperbolic conservation laws with greater efficiency than first- and second-order discretizations.  However, high-order methods tend to be less robust, and this has hindered their widespread adoption.  The need for robust and efficient high-order discretizations motivates this work on entropy-stable discontinuous Galerkin difference methods.  

Our interest in the Galerkin difference (GD) family of methods stems from the attractive properties of their underlying basis functions. For example, in the original GD method proposed by Banks and Hagstrom~\cite{Banks:GD2016}, the solution is expressed in terms of piecewise continuous functions that extend over multiple elements.  Consequently, the number of GD degrees of freedom remains constant as the polynomial degree increases, similar to conventional finite-difference and finite-volume methods.  Furthermore, like finite-difference methods, GD schemes have time-step restrictions that are relatively insensitive to their order of accuracy.

The GD method was recently generalized to discontinuous basis functions~\cite{Banks-DGD-paper}, and this discontinuous Galkerin difference (DGD) method is the starting point for the discretization considered in this work.  Building on~\cite{Banks-DGD-paper}, we present an entropy-stable formulation of DGD to address robustness for discretizations of symmetrizable hyperbolic systems.  In addition, unlike~\cite{Banks-DGD-paper}, we consider unstructured grids and construct basis function stencils following the approach presented by Li~\etal~\cite{Li2019}.

To derive the entropy-stable DGD method, we leverage the existing theory for entropy-stable summation-by-parts (SBP) methods.  Over the past decade, researchers have used SBP operators to construct high-order entropy-stable finite-difference~\cite{fisher:thesis,fisher:2013,Fisher2013discretely}, finite-element~\cite{Chan2018discretely}, and spectral-element-type methods~\cite{Carpenter2014entropy,parsani:2016,Gassner2016well,Ranocha2018stability,Friedrich2018entropy,Friedrich2019entropy,Shadpey2020entropy,Rojas2021robustness}.  Entropy-stable discretizations provide a form of nonlinear stability and, therefore, robustness by ensuring the entropy is non-increasing\footnote{Entropy here refers to mathematical entropy.}.

In summary, the objectives of this work are to present a framework for constructing entropy-stable DGD discretizations and to study the properties of the resulting discretization.  The principle contributions are listed below:
\begin{itemize}
    \item DGD discretizations of linear and nonlinear hyperbolic systems can be constructed using diagonal-norm SBP operators.
    \item DGD matrix operators are, themselves, dense-norm SBP operators.
    \item When applied to entropy-stable DGD semi-discretizations, explicit time marching schemes require the solution of a coupled nonlinear algebraic system and, therefore, do not offer a significant computational advantage over implicit methods.
\end{itemize}

The roadmap of the paper is as follows. Section \ref{GDmethod} reviews the DGD method, including the construction of the DGD basis functions and the semi-discretization of the two-dimensional linear advection problem.  Section \ref{GDSBP-relation} reviews multi-dimensional SBP operators and shows how they can be used in the DGD semi-discretization of linear advection.  This section also proves that DGD operators are dense-norm SBP operators.  In Section \ref{sec:entropytheory}, we review generic entropy-conservative/stable SBP discretizations, and then develop the entropy-stable DGD discretization.  Section \ref{numericalexp} presents numerical experiments in order to verify the accuracy and stability properties of the DGD discretization, as well as characterize its efficiency.  Finally, Section \ref{conclusion} concludes this study with a summary and a discussion of potential future developments.

\section{The discontinuous Galerkin difference method on unstructured grids}
\label{GDmethod}
Galerkin difference methods are finite element methods that use piecewise polynomial basis functions.  One of our objectives in this section is to familiarize readers with these somewhat unconventional basis functions.  We also highlight some important differences that arise when DGD basis functions are constructed on unstructured versus tensor-product grids and when the polynomial degree is not even.  We conclude this section by illustrating how the DGD basis functions are used to discretize the constant-coefficient linear advection equation.
\subsection{DGD basis functions}
\label{sec:basis}
Let $T_h = \{ \Omega_k \}_{k=1}^K$ denote a tesselation of a closed and bounded domain $\Omega$ into $K$ elements, in which element $k$ has subdomain $\Omega_k$, boundary $\Gamma_k$, and centroid $\tilde{\bm{x}}_k$. The stencil $N_k = \{ \nu_1, \nu_2, \ldots, \nu_{n_k} \} \subset \{1,2,3, ... K \}$ represents the ordered set of elements whose basis functions are nonzero on element $k$ and, therefore, influence the discrete solution over $\Omega_k$.  We will sometimes refer to the stencil $N_k$ as the patch for element $k$.  The number of elements in the stencil is denoted by $n_k$.

\begin{figure}[tbp]
\begin{center}
  \subfigure[$p=1$]{%
        \includegraphics[width=0.46\textwidth]{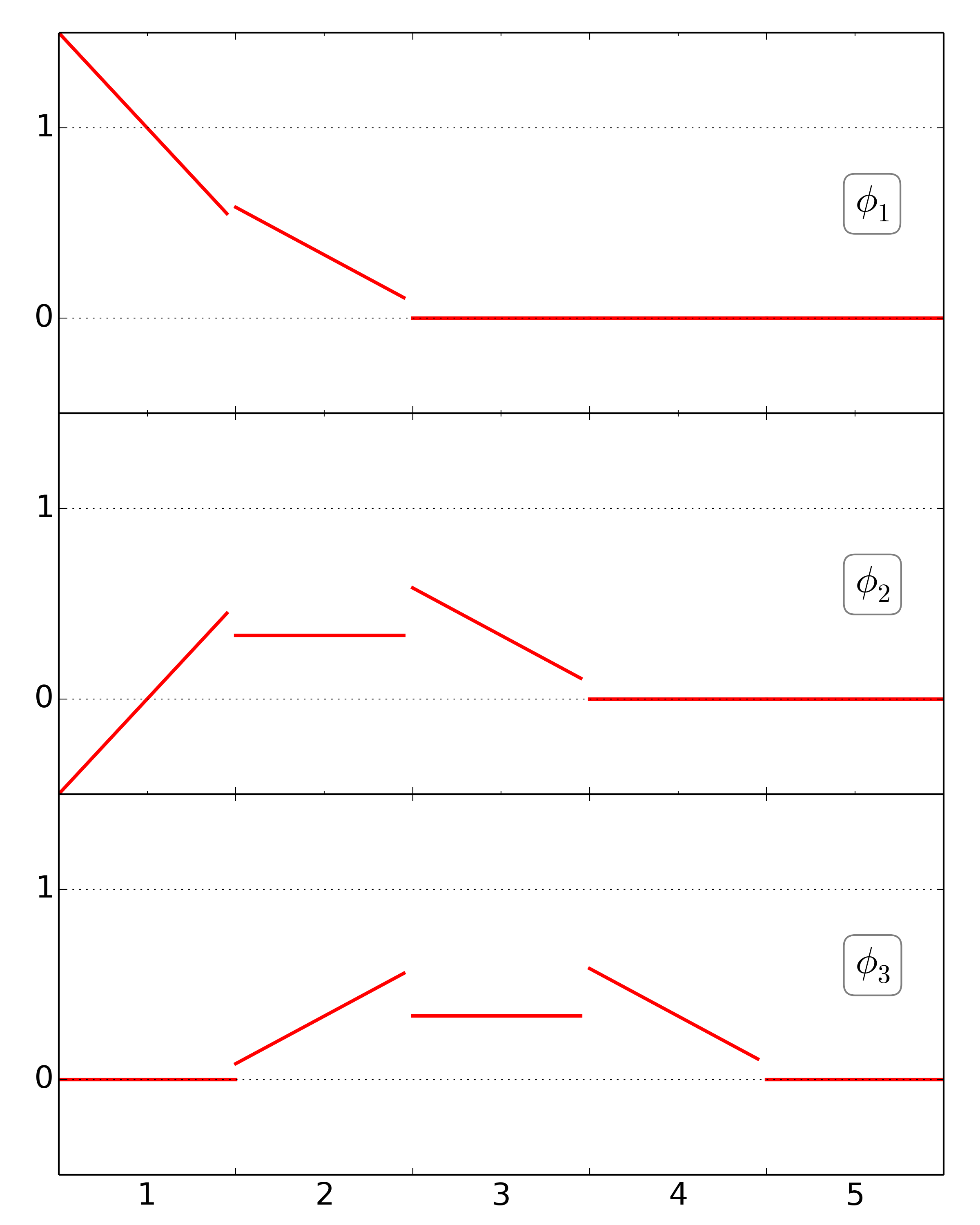}}\hfill
  \subfigure[$p=2$]{%
        \includegraphics[width=0.46\textwidth]{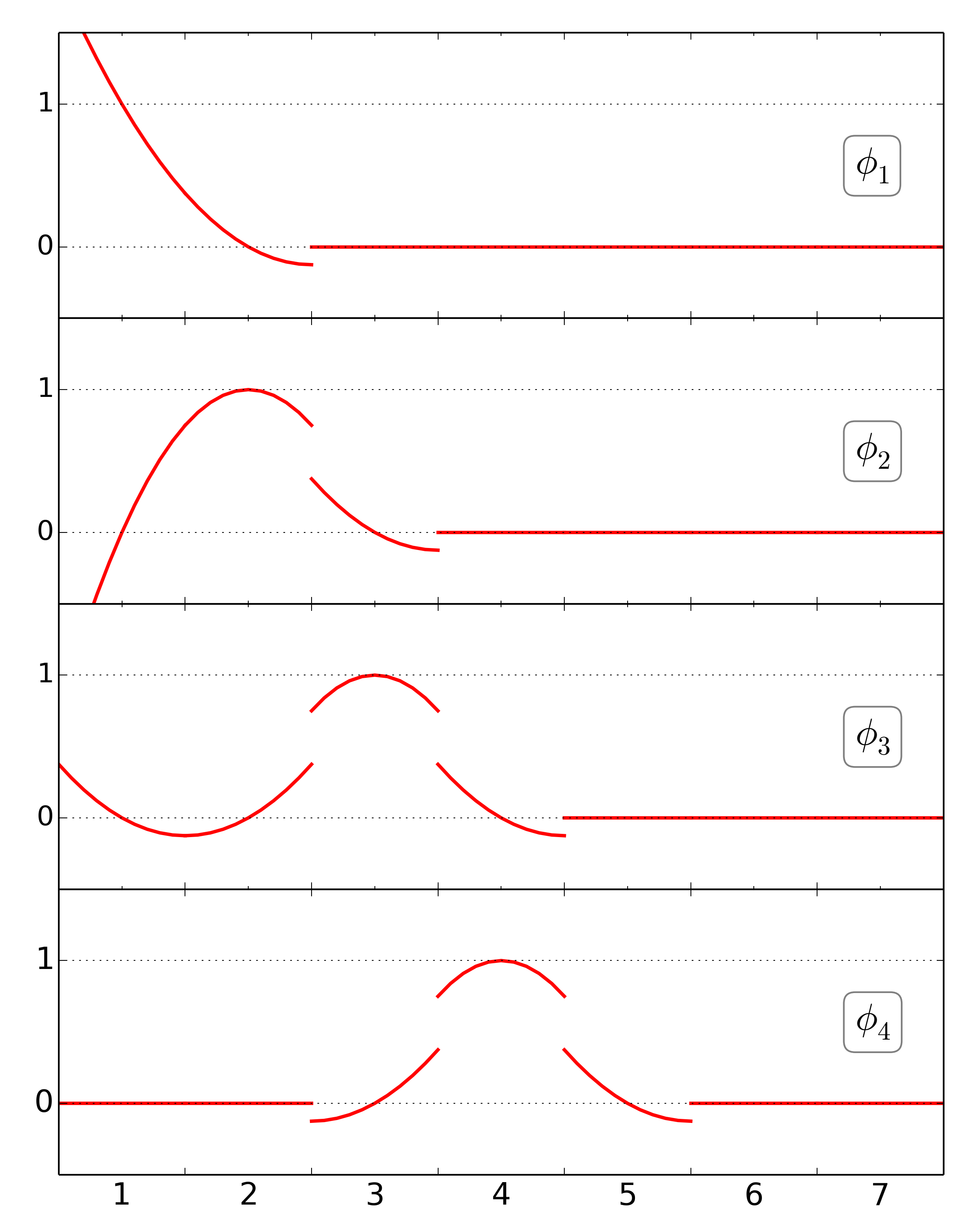}} \\
  \subfigure[$p=3$]{%
        \includegraphics[width=0.46\textwidth]{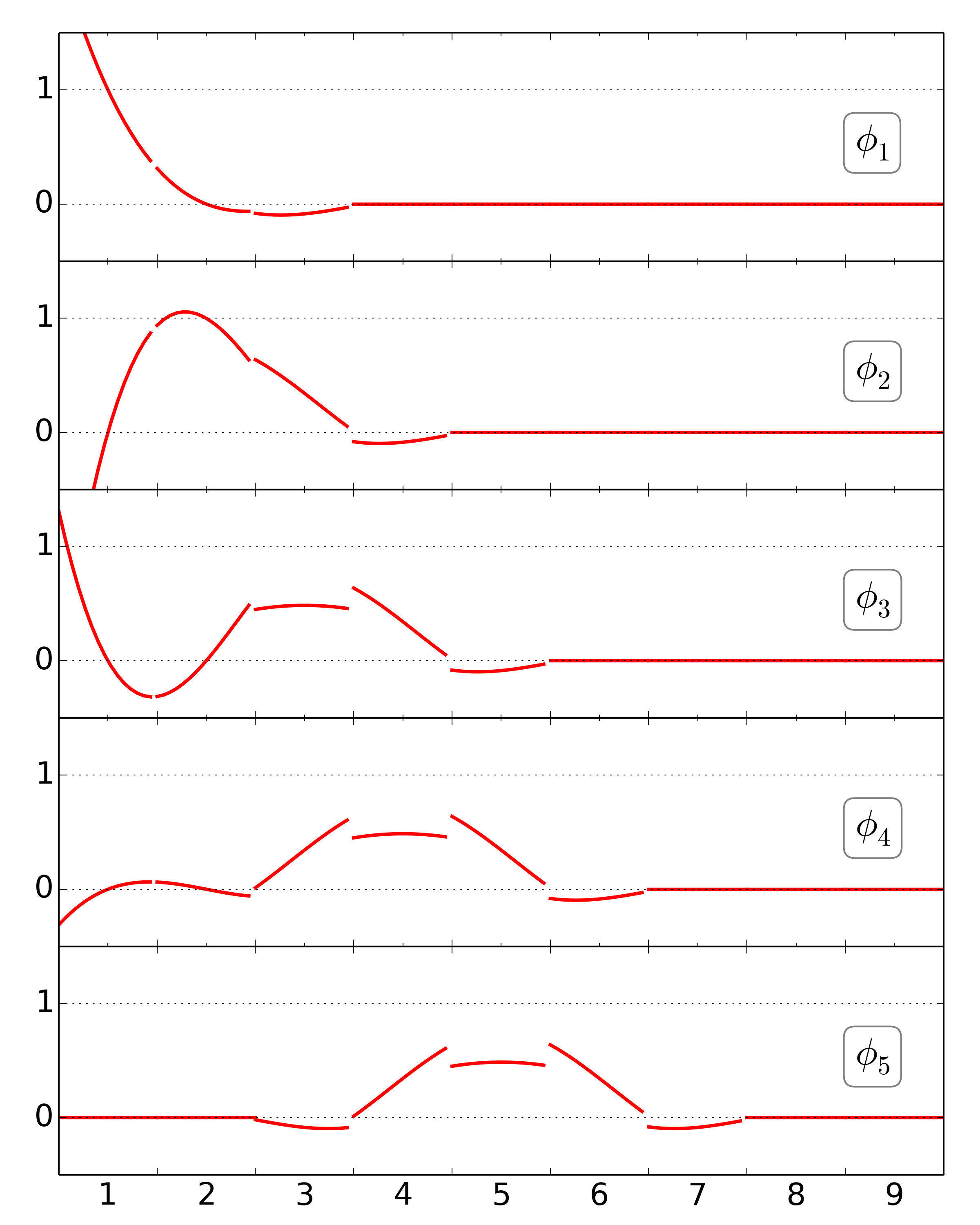}}\hfill
  \subfigure[$p=4$]{%
        \includegraphics[width=0.46\textwidth]{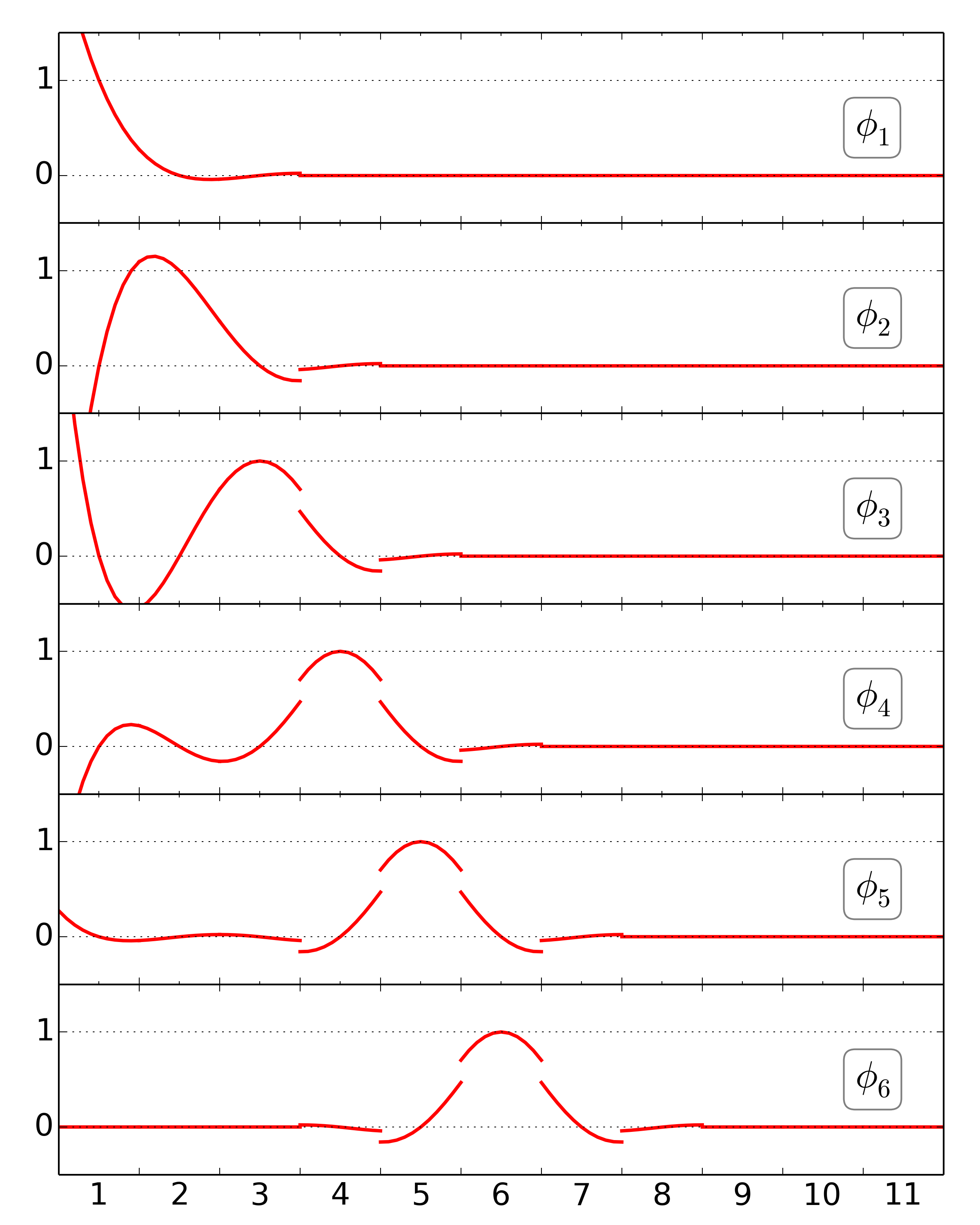}}
\end{center}
  \caption{Discontinuous Galerkin difference basis functions, for degrees $p=1$, $p=2$, $p=3$, and $p=4$, on a uniform one dimensional grid.  The $x$-axis labels indicate the element indices.  Only the basis functions on the left half of the domain are shown, since the right-half basis functions are symmetric. \label{fig:dgd_basis}}
\end{figure}

Consider a bounded function $\mathcal{U}(\bm{x})$ and let $\tilde{u}_\nu = \mathcal{U}(\tilde{\bm{x}}_\nu)$ denote its value at the centroid of element $\nu$.  Then the DGD interpolation of $\mathcal{U}$, which we denote by $\mathcal{U}_h(\bm{x})$, takes the following form on element $k$:
\begin{equation}\label{eq:localinterp}
    \mathcal{U}_h(\bm{x}) = 
    \sum_{\nu=1}^{K} \tilde{u}_\nu \phi_\nu(\bm{x}) = 
    \sum_{\nu \in N_k} \tilde{u}_{\nu} \phi_{\nu}(\bm{x}), \quad \bm{x} \in \Omega_k,
\end{equation}
where $\phi_\nu(\bm{x})$ is the piecewise polynomial basis associated with element $\nu$ whose construction is the focus of this section.  To help readers gain some intuition about DGD basis functions, some one dimensional examples of $\phi_\nu(\bm{x})$ are illustrated in  Figure~\ref{fig:dgd_basis}.

In general, the DGD basis functions can be expressed as a linear combination of some (standard) polynomial basis functions on each element.  Specifically, the basis of the $i$th neighbour in the stencil $N_k$ can be written as
\begin{equation}\label{eq:basis_general}
\phi_{\nu_i}(\bm{x}) = \sum_{j = 1}^{n_p} c_{ji} \mathcal{V}_j(\bm{x}),
\qquad \forall \bm{x} \in \Omega_k,
\end{equation}
where $\{\mathcal{V}_j\}_{j=1}^{n_p}$ is a given basis for $\mathbb{P}_p(\Omega_k)$, the space of total degree $p$ polynomials on $\Omega_k$, and $n_p$ is the dimension of the basis.

Ideally, the coefficients $c_{ji}$ are chosen such that the DGD basis $\phi_{\nu_i}(\bm{x})$ has a value of one at the centroid $\tilde{\bm{x}}_{\nu_i}$ and a value of zero at the centroids of the other elements in $N_k$.  In other words, for all $\nu_i, \nu_m \in N_k$, we have
\begin{equation} \label{eq:ploybasis}
    \phi_{\nu_i}(\Tilde{\bm{x}}_{\nu_m}) = \sum_{j=1}^{n_p} c_{ji}\mathcal{V}_j(\Tilde{\bm{x}}_{\nu_m}) = \delta_{im},
\end{equation}
where $\delta_{im}$ denotes the Kronecker delta.  Equation \eqref{eq:ploybasis} can be expressed in matrix form as
\begin{equation}\label{eq:basis_coeff}
    \Vtilde \mat{C}_k = \mat{I}_k,
\end{equation}
where $\mat{I}_k$ is the $n_k \times n_k$ identity matrix and $\mat{C}_k \in \mathbb{R}^{n_p \times n_k}$ holds the to-be-determined coefficients. $\Vtilde \in \mathbb{R}^{n_k \times n_p}$ is the generalized Vandermonde matrix:
\begin{equation*}
\Vtilde \equiv \begin{bmatrix}
  \mathcal{V}_{1}(\tilde{\bm{x}}_{\nu_1}) & \mathcal{V}_{2}(\tilde{\bm{x}}_{\nu_1}) & \hdots 
  & \mathcal{V}_{n_p}(\tilde{\bm{x}}_{\nu_1}) \\
  \mathcal{V}_{1}(\tilde{\bm{x}}_{\nu_2}) & \mathcal{V}_{2}(\tilde{\bm{x}}_{\nu_2}) & \hdots 
  & \mathcal{V}_{n_p}(\tilde{\bm{x}}_{\nu_2}) \\
  \vdots & \vdots & \ddots & \vdots \\
  \mathcal{V}_{1}(\tilde{\bm{x}}_{\nu_{n_k}}) & \mathcal{V}_{2}(\tilde{\bm{x}}_{\nu_{n_k}}) & \hdots 
  & \mathcal{V}_{n_p}(\tilde{\bm{x}}_{\nu_{n_k}})
\end{bmatrix}.
\end{equation*}

In order for \eqref{eq:basis_coeff} to have a unique solution for $\mat{C}_k$, the number of elements in the stencil must be the same as the number of basis functions, $n_k = n_p$, and the element centroids must form an unisolvent point set for $\{\mathcal{V}_j\}_{j=1}^{n_p}$. Assuming the degrees of freedom are stored at element centers (\ie, on the dual grid), the constraint $n_k = n_p$ is straightforward to enforce for tensor-product DGD elements of even degree $p$, but it is problematic for odd $p$ and unstructured grids.  For instance, only three elements are required in the stencil to ensure $n_k = n_p$ when constructing a piecewise linear ($p=1$) total degree DGD basis in two dimensions; however, on a triangular mesh, restricting $|N_k|$ to three would require excluding either the element $k$ itself or one of its adjacent neighbors.  As $p$ increases, the choice of which elements to exclude from the stencil $N_k$ becomes increasingly arbitrary.

In this work, we favor increasing $n_k$ rather than excluding particular neighbors and potentially biasing the stencil.  One consequence of this is that we sacrifice interpolation when $n_k > n_p$, since we solve \eqref{eq:basis_coeff} for $\mat{C}_k$ in a least squares sense when there are more elements in the stencil than basis functions:
\begin{equation} \label{eq:csolution}
    \mat{C}_k = (\Vtilde^{T} \Vtilde)^{-1}\Vtilde^{T}.
\end{equation}
Thus, unlike the tensor-product scheme in~\cite{Banks-DGD-paper}, the present DGD solution does not necessarily interpolate data at the element centroids; nevertheless, the basis is still capable of representing degree $p$ polynomials exactly. To prove this (see also \cite{Li2019}), it is sufficient to show that the DGD basis can represent any basis function $\fnc{V}_{m} \in \{\mathcal{V}_j\}_{j=1}^{n_p}$.  Indeed, if we choose $\tilde{u}_{\nu_i} = [\Vtilde]_{im} = \fnc{V}_{m}(\tilde{\bm{x}}_{\nu_i})$ in the DGD-basis expansion \eqref{eq:localinterp} then, for all $\bm{x} \in \Omega_k$, we find
\begin{equation}
\begin{split}
    \fnc{U}_h(\bm{x}) = \sum_{\nu \in N_k} \tilde{u}_\nu \phi_\nu(\bm{x}) &= \sum_{i=1}^{n_k} [\Vtilde]_{im} \sum_{j=1}^{n_p} c_{ji} \mathcal{V}_j(\bm{x}) \\ 
    &= \sum_{j=1}^{n_p} \mathcal{V}_j(\bm{x}) \underbrace{\sum_{i=1}^{n_k} [(\Vtilde^T \Vtilde)^{-1} \Vtilde^T ]_{ji} [\Vtilde]_{im}}_{\displaystyle = \delta_{jm}}
    = \fnc{V}_{m}(\bm{x}).
\end{split}
\end{equation}
It follows that the DGD basis can represent any polynomial of total degree $p$ on the entire domain $\Omega$, since the element $k$ is arbitrary.

One-dimensional DGD basis functions of degree $p=1$ to $p=4$ are illustrated in Figure \ref{fig:dgd_basis}. Notice that the even degrees $p=2$ and $p=4$ produce interpolatory basis functions, with $\phi_\nu(\bm{x})$ equal to one at $\tilde{\bm{x}}_\nu$ and zero at the other nodes.  This is because $n_k = n_p$ for these basis functions, so \eqref{eq:basis_coeff} has a unqiue solution.  In contrast, the basis functions $\phi_\nu(\bm{x})$ for the odd degrees do not equal one at $\tilde{\bm{x}}_\nu$ since a symmetric stencil is enforced on the interior elements, resulting in $n_k > n_p$.

\subsection{Reconstruction operator}

In practice, we only need to evaluate the DGD basis at quadrature points in order to compute the integrals that arise in the finite-element method.  Thus, in this section, we show how the basis functions can be used to construct a linear mapping from the centroids in $N_k$ to quadrature points.  We will refer to this linear mapping as the prolongation matrix; in the finite-volume literature this matrix is often called the reconstruction operator.

Let $X_k \equiv \{ \bm{x}_q\}_{q = 1}^{n_q}$ denote the quadrature points adopted for the element domain $\Omega_k$. Then the discrete solution at $\bm{x}_q \in X_k$ is given by
\begin{equation}\label{eq:s22localinterp}
    \mathcal{U}_h(\bm{x}_q) = \sum_{\nu \in N_k} \phi_{\nu}(\bm{x}_q) \tilde{u}_{\nu}
    = \sum_{\nu=1}^{K} (\mat{P}_k)_{q\nu} \tilde{u}_{\nu}, \quad \forall q = 1,2, \dots , n_q,
\end{equation}
where $(\mat{P}_k)_{q\nu} \equiv \phi_\nu(\bm{x}_q)$ is simply the basis function corresponding to element $\nu$ evaluated at the quadrature point $\bm{x}_q$.  Equation~\eqref{eq:s22localinterp} can be written succinctly as 
\begin{equation}
\mat{P}_k \ugd = \bm{u}_k,
\end{equation}
where $\mat{P}_k \in \mathbb{R}^{n_q \times K}$ is the prolongation matrix for element $k$, $\ugd \in \mathbb{R}^{K}$ is a vector that holds DGD coefficients associated with the centroids of the elements, and $\bm{u}_k$ is $\mathcal{U}_h$ evaluated at, or ``prolonged'' to, the quadrature points.

To avoid constructing the DGD basis explicitly, we express the prolongation matrix $\mat{P}_k$ in terms of the polynomial basis $\{\fnc{V}_j\}_{j=1}^{n_p}$.  The desired expression can be inferred from \eqref{eq:s22localinterp}, \eqref{eq:csolution}, and \eqref{eq:basis_general}:
\begin{equation}\label{eq:psolution}
    \mat{P}_k = \mat{V}_k \left( \Vtilde^T \Vtilde \right)^{-1} \Vtilde^{T} \mat{Z}_k,
\end{equation}
where $\mat{V}_k \in \mathbb{R}^{n_q \times n_p}$ and its entries $[\mat{V}_k]_{qj} = \mathcal{V}_j(\bm{x}_q)$ are the polynomial basis evaluated at the quadrature points $X_k$.  Furthermore, the matrix $\mat{Z}_k \in \mathbb{R}^{n_k \times K}$ maps from the global indices $\{1,2,\ldots,K\}$ to the local indices $\{1,2,\ldots,n_k\}$ of $N_k$:
\begin{equation*}
    (\mat{Z}_{k})_{i\nu} = \begin{cases}
      1, & \nu = \nu_i \in N_k, \\
      0, & \text{otherwise}.
    \end{cases}
\end{equation*}

\begin{remark}\label{rmk:prolong_direct}
The prolongation matrix can be constructed directly by seeking an interpolation operator from the centroids of the stencil $N_k$ to the quadrature points that is exact for all degree $p$ polynomials. The conditions for such an operator are given by
\begin{equation}\label{eq:prolong_direct}
\mat{P}_k \mat{Z}_k^T \Vtilde = \mat{V}_k.
\end{equation}
This equation is underdetermined when $n_q \geq n_p$ --- assuming a unisolvent set of quadrature points --- so additional conditions are necessary to fix a unique $\mat{P}_k$.  Here we seek the prolongation matrix that has the minimum Frobenius norm and satisfies \eqref{eq:prolong_direct}.  It is straightforward to show that such a matrix is given by \eqref{eq:psolution}.
\end{remark}

\subsection{DGD discretization of the linear advection equation}
\label{sec:DGD_linear_advection}

We conclude this section by applying the DGD method to semi-discretize the two-dimensional linear advection equation.  This exposition is intended to further familiarize readers with the DGD finite-element method, but it will also be used later to relate the method to summation-by-parts discetizations and construct entropy-conservative/stable DGD schemes.

Consider the two dimensional, constant-coefficient linear advection equation on the domain $\Omega$:
\begin{equation}\label{eq:simpleeg}
        \frac{\partial \mathcal{U}}{\partial t} + \lambda_x \frac{\partial \mathcal{U}}{\partial x} + \lambda_y \frac{\partial \mathcal{U}}{\partial y} = 0,\qquad \forall\; (x,y) \in \Omega,
\end{equation}
where $\bm{\lambda} = [\lambda_x,\lambda_y]^T$ is the advection velocity.  In practice, the PDE \eqref{eq:simpleeg} requires boundary conditions and an initial condition; however, given our present focus on the spatial operators, we will ignore the boundary and initial conditions for the time being.


The DGD semi-discretization of the linear advection equation is obtained by following the usual weighted-residual approach.  Let $W_h \equiv \spn \{ \phi_k \}_{k=1}^K$ denote the DGD finite-dimensional function space. 
Then the DGD weak formulation seeks $\mathcal{U}_h \in W_h$ such that
\begin{multline}\label{eq:weakform}
    \sum_{k=1}^{K} \Bigg[ \int_{\Omega_k} \mathcal{V}_h \frac{\partial \mathcal{U}_h}{\partial t} \, d\Omega 
    - \int_{\Omega_k} \lambda_x \frac{\partial \mathcal{V}_h}{\partial x} \mathcal{U}_h \, d\Omega 
    - \int_{\Omega_k} \lambda_y \frac{\partial \mathcal{V}_h}{\partial y} \mathcal{U}_h \, d\Omega \\
    + \int_{\Gamma_k} \mathcal{V}_h \mathcal{U}_h \lambda_x n_x \, d\Gamma 
    + \int_{\Gamma_k} \mathcal{V}_h \mathcal{U}_h \lambda_y n_y \, d\Gamma \Bigg] = 0,
\end{multline}
for all $\displaystyle \mathcal{V}_h \in W_h$.  The integrals are taken over the elements and their boundaries to accommodate the discontinuous basis functions.  


As usual with finite-element methods, the bilinear forms in the weak formulation \eqref{eq:weakform} can be represented as matrices.  To find this equivalent representation, we express the trial and test functions as
\begin{equation*}
    \mathcal{U}_h(\bm{x}) = \sum_{k=1}^{K} \tilde{u}_k \phi_k(\bm{x}),
    \qquad\text{and}\qquad
    \mathcal{V}_h(\bm{x}) = \sum_{k=1}^{K} \tilde{v}_k \phi_k(\bm{x})
\end{equation*}
and substitute these expansions into \eqref{eq:weakform}.  Thus, the DGD weak formulation is equivalent to 
\begin{equation}\label{eq:dis_weak_GD}
  \vgd^T \tM \frac{d \ugd}{dt}
  - \lambda_x \vgd^T \tQx^T \ugd - \lambda_y \vgd^T \tQy^T \ugd  
  + \lambda_x \vgd^T \tEx \ugd + \lambda_y \vgd^T \tEy \ugd = 0,
\end{equation}
for all $\vgd \in \mathbb{R}^{K}$, 
where $\ugd \in \mathbb{R}^{K}$ is the vector of solution coefficients, and the matrices $\tM$, $\tQx$, and $\tEx$ are defined by
\begin{equation}\label{entries}
\begin{gathered}
\tM_{ij} = \sum_{k=1}^{K} \int_{\Omega_k} \phi_i \phi_{j}\, d\Omega,
\qquad
(\tQx)_{ij} = \sum_{k=1}^{K} \int_{\Omega_k} \phi_i \frac{\partial \phi_j}{\partial x} \, d\Omega, \\ 
\text{and}\qquad
(\tEx)_{ij} = \sum_{k=1}^{K} \int_{\Gamma_k} \phi_i \phi_j n_x  \, d\Gamma.
\end{gathered}
\end{equation}
The matrices $\tQy$ and $\tEy$ are defined analogously to $\tQx$ and $\tEx$, respectively.

\begin{remark}
    The DGD mass matrix $\tM$ is sparse, since $\tM_{ij}$ is non-zero only if both elements $i$ and $j$ are included in a common stencil: that is, there exists $k \in \{1,2,\ldots,K\}$ such that $i,j \in N_k$.  For the same reason, the operators $\tQx$, $\tQy$, $\tEx$, and $\tEy$ are also sparse.
\end{remark}

\begin{remark}    
    Unlike DG mass matrices, the DGD mass matrix is not block diagonal, which has implications for explicit time marching methods.  However, on structured grids, Galerkin difference methods can take advantage of the tensor-product structure to invert the mass matrix rapidly in linear time \cite{banks:GD_highorder2019}.
\end{remark}

\section{DGD and summation-by-parts operators}
\label{GDSBP-relation}

The goal of this section is two-fold.  First, we will show how DGD discretizations can be implemented using summation-by-parts (SBP) operators and the element prolongation matrices, $\mat{P}_k$.  Second, we will prove that DGD operators are themselves dense-norm SBP operators.

\subsection{Multidimensional SBP operators and their properties}

Consider one of the elements $\Omega_k \in T_h$, and let $X_k = \{ \bm{x}_q \}_{q=1}^{n_q}$ be a set of nodes that are in the closure of the element subdomain, $\bm{x}_q \in \bar{\Omega}_k$, $\forall\, q=1,2,\ldots,n_q$.  The notation for the nodes $X_k$ is the same used earlier for quadrature points; as we shall see, this choice is deliberate.

\begin{definition}[Summation-by-parts first-derivative operator]
The matrix $\Dxk \in \mathbb{R}^{n_q \times n_q}$ is a degree $p$ summation-by-parts operator approximating the first-derivative with respect to $x$ at the nodes $X_k$ if the following conditions are satisfied.
\begin{enumerate}
    \item The difference operator exactly differentiates polynomials of total degree $p$ at the nodes $X_k$:
    \begin{equation}\label{eq:SBP_accuracy}
    \sum_{q=1}^{n_q} (\Dxk)_{rq} \fnc{V}_i(\bm{x}_q) = \frac{\partial \fnc{V}_i}{\partial x}(\bm{x}_r),
    \qquad\forall\, \bm{x}_r \in X_k,
    \end{equation}
    and for all polynomials $\fnc{V}_i$ in the basis $\{ \fnc{V}_j \}_{j=1}^{n_p} \subset \mathbb{P}_p(\Omega_k)$.
    \item $\Dxk = \Hk^{-1} \Qxk$, where $\Hk$ is a symmetric positive-definite matrix.
    \item $\Qxk = \Sxk + \frac{1}{2} \Exk$, where $(\Sxk)^T = -\Sxk$ is a skew-symmetric matrix, and the symmetric matrix $\Exk = (\Exk)^T$ satisfies
    \begin{equation}\label{eq:SBP_Ex_accuracy}
        \sum_{r=1}^{n_q} \sum_{q=1}^{n_q} \fnc{V}_i(\bm{x}_r) (\Exk)_{rq} \fnc{V}_j(\bm{x}_q) = \int_{\Gamma_k} \fnc{V}_i \fnc{V}_j n_x \, d\Gamma,
    \end{equation}
    for all basis polynomials $\fnc{V}_i, \fnc{V}_j \in \{ \fnc{V}_j \}_{j=1}^{n_p} \subset \mathbb{P}_p(\Omega_k)$, where $n_x$ is the $x$ component of the outward pointing unit normal on $\Gamma_k$.
\end{enumerate}
An analogous definition holds for the SBP difference operator $\Dyk$, which approximates the first-derivative in the $y$ direction. 
\end{definition}

An important subset of SBP operators, called diagonal-norm operators, are worth highlighting and will be used later.  Diagonal-norm operators have diagonal $\Hk$, and one can show that the diagonal entries $\{ (\Hk)_{qq} \}_{q=1}^{n_q}$ and nodes $X_k$ constitute a quadrature rule that is at least $2p-1$ exact~\cite{hicken:quad2013,hicken:mdimsbp2016}.  That is, we have
\begin{equation*}
    \sum_{q=1}^{n_q} (\Hk)_{qq} \fnc{V}_i(\bm{x}_q) = \int_{\Omega_k} \fnc{V}_i \, d\Omega,
\end{equation*}
for all $\fnc{V}_i \in \{ \fnc{V}_j \}_{j=1}^{n_{2p-1}} \subset \mathbb{P}_{2p-1}(\Omega_k).$  Note that the accuracy of $2p-1$ exactness is a lower bound.  For the subsequent analysis, we will assume that we are using diagonal-norm SBP operators whose quadrature accuracy is at least $2p$ exact.

\subsection{Implementation of DGD with SBP operators}

In this section, we use diagonal-norm SBP operators to discretize the constant-coefficient linear-advection equation, and then show the relationship between this SBP discretization and the DGD discretization.

Suppose we have SBP operators $\Dxk$ and $\Dyk$ for each element $k \in \{ 1,2,\ldots,K \}$ in the tesselation $T_h$ of $\Omega$.  Then, the SBP discretization of the linear-advection equation~\eqref{eq:simpleeg} is given by 
\begin{equation}\label{eq:SBP_linear_advection_strong}
    \frac{d \uk}{d t} + \lambda_x \Dxk \uk + \lambda_y \Dyk \uk = \bm{0},\qquad \forall \, k = 1,2,\ldots,K,
\end{equation}
where $\uk \in \mathbb{R}^{n_q}$ denotes the SBP solution at the nodes $X_k$ of element $k$.  As with the DGD semi-discretization in Section~\ref{sec:DGD_linear_advection}, this SBP discretization is for illustrative purposes only, since it lacks imposition of boundary conditions and inter-element coupling.

Equation \eqref{eq:SBP_linear_advection_strong} is the strong form of the SBP discretization.  The equivalent weak form can be obtained by left multiplying by $\vk^T \Hk$, where $\vk \in \mathbb{R}^{n_q}$ denotes a test function at the nodes $X_k$, and summing over all elements:
\begin{multline}\label{eq:SBP_linear_advection_weak}
    \sum_{k=1}^{K} \Bigg[ \vk^T \Hk \frac{d \uk}{d t} - \lambda_x \vk^T \Qxk^T \uk - \lambda_y \vk^T \Qyk^T \uk \\
    + \lambda_x \vk^T \Exk \uk + \lambda_y \vk^T \Eyk \uk \Bigg] = 0,
\end{multline}
for all $\vk \in \mathbb{R}^{n_q}$ and all 
$k = 1,2,\ldots,K$.  To arrive at  \eqref{eq:SBP_linear_advection_weak}, we used the identities $\Qxk + \Qxk^T = \Exk$ and $\Qyk + \Qyk^T = \Eyk$.

Our goal is to relate the SBP discretization \eqref{eq:SBP_linear_advection_weak} and the DGD discretization \eqref{eq:dis_weak_GD}.  To that end, we will need the following lemma, which provides identities relating the SBP and DGD matrix operators.

\begin{lemma}\label{lem:DGD_SBP_matrices}
For each element $k \in \{ 1,2,\ldots,K \}$, let $\Dxk$ be a diagonal-norm SBP operator whose corresponding quadrature rule is at least $2p$ exact.  In addition, let $\mat{P}_k$ be the DGD prolongation matrix defined in \eqref{eq:s22localinterp}.
Then the DGD matrices in \eqref{entries} can be evaluated as
\begin{equation}\label{eq:DGD_SBP_identities}
    \begin{gathered}
    \tM_{ij} = \left[ \sum_{k=1}^{K} \mat{P}_k^T \Hk \mat{P}_k \right]_{ij},
    \qquad 
    (\tQx)_{ij} = \left[ \sum_{k=1}^K \mat{P}_k^T \Qxk \mat{P}_{k} \right]_{ij}, \\
    \text{and}\qquad 
    (\tEx)_{ij} = \left[ \sum_{k=1}^K \mat{P}_k^T \Exk \mat{P}_{k} \right]_{ij}.
    \end{gathered}
\end{equation}
\end{lemma}

The proof of Lemma~\ref{lem:DGD_SBP_matrices} can be found in~\ref{app:DGD_SBP_proof}.  Note that an analogous result exists for the $y$-coordinate operators, $\tQy$ and $\tEy$.

We can now state and prove our first result.

\begin{thrm}\label{thm:DGD_SBP_equiv}
Let $\Dxk$ and $\Dyk$ be diagonal-norm SBP operators for element $k \in \{ 1,2,\ldots,K \}$, and assume the quadrature associated with these SBP operators is $2p$ exact.  If we define $\uk = \mat{P}_k \ugd$ and $\vk = \mat{P}_k \vgd$ on each element, where $\mat{P}_k$ is the DGD prolongation matrix defined in \eqref{eq:s22localinterp}, the SBP discretization~\eqref{eq:SBP_linear_advection_weak} of the constant-coefficient linear advection is equivalent to the DGD discretization~\eqref{eq:dis_weak_GD}.
\end{thrm}

\begin{proof}
The result will follow if we can show that each term in \eqref{eq:SBP_linear_advection_weak} is equivalent to the corresponding term in \eqref{eq:dis_weak_GD}.  We will establish this equivalence for the temporal term.

Substituting $\uk = \mat{P}_k \ugd$ and $\vk = \mat{P}_k \vgd$ into the temporal term we find
\begin{equation*}
\sum_{k=1}^{K}
\vk^T \Hk \frac{d \uk}{d t}
= \sum_{k=1}^{K} \vgd^T \mat{P}_k^T \Hk \frac{d}{dt} \left( \mat{P}_k \ugd \right)
= \vgd^T \left[ \sum_{k=1}^{K} \mat{P}_k^T \Hk \mat{P}_k \right] \frac{d \ugd}{dt}
= \vgd^T \tM \ugd, 
\end{equation*}
where we used Lemma~\ref{lem:DGD_SBP_matrices} in the final step.  This demonstrates that the temporal terms are equivalent under the assumptions.  The result follows after applying a similar substitution in the remaining terms and applying Lemma~\ref{lem:DGD_SBP_matrices}.
\end{proof}

Theorem~\ref{thm:DGD_SBP_equiv} is useful because it provides a procedure for constructing DGD discretizations from SBP discretizations.  This relationship between the methods can be exploited for proofs of entropy conservation and stability, as we will show in Section~\ref{sec:entropytheory}.

\subsection{DGD operators are dense-norm SBP operators}

The previous section described how DGD discretizations can be implemented using SBP operators.  Here we explore a closely related connection between these discretizations; namely, that DGD difference operators are, themselves, a type of SBP operator.

\begin{thrm}\label{thm:DGD_operators}
Let $T_h = \{\Omega_k \}_{k=1}^K$ denote a tesselation of the domain $\Omega$ into non-overlapping elements.  Assume that a degree $2p$ exact quadrature rule with positive weights and nodes $X_k$ is available on each element $k$, and that 
\begin{equation*}
    \rank\left(\left[\mat{P}_1^T,\mat{P}_2^T,\ldots,\mat{P}_K^T\right]\right) = K,
\end{equation*}
where $\mat{P}_k$ is the prolongation matrix from the element centers to $X_k$.  Then the degree $p$ DGD operator $\tDx = \tM^{-1} \tQx$ is a degree $p$ multidimensional SBP operator for the partial derivative in the $x$ direction over the domain $\Omega$.
\end{thrm}

\begin{proof}
The proof relies on Lemma~\ref{lem:DGD_SBP_matrices} and, therefore, the availability of degree $p$ diagonal-norm SBP operators, $\Dxk$, with quadrature accuracy at least $2p$ on each element.  The existence of such operators is guaranteed by the assumption that there is a sufficiently accurate quadrature rule with positive weights for each element~\cite{hicken:mdimsbp2016}.

Let $\vgd_i \in \mathbb{R}^K$ denote an arbitrary basis polynomial $\fnc{V}_i \in \{ \fnc{V}_j \}_{j=1}^{n_p} \subset \mathbb{P}_{p}(\Omega)$, evaluated at the centroids of the elements in $T_h$, and $\vgd_i' \in \mathbb{R}^K$ denote its derivative, $\partial \fnc{V}_i/\partial x$, also evaluated at the centroids of the elements.  Similarly, let $\bm{v}_{i,k}$ and $\bm{v}_{i,k}'$ denote the basis function and its derivative, respectively, evaluated at the SBP nodes $X_k$ of element $k$.

We begin with the SBP accuracy property~\eqref{eq:SBP_accuracy}, which, after multiplying both sides of the identity by $\tM$, is equivalent to
\begin{equation*}
    \tQx \vgd_i = \tM \vgd_i',
\end{equation*}
Consider the left-hand side of the above equation.  Substituting the expression for $\tQx$ from \eqref{eq:DGD_SBP_identities}, we find 
\begin{align*}
    \tQx \vgd_i &=
    \left[ \sum_{k=1}^{K} \mat{P}_k^T \Qxk \mat{P}_k \right] \vgd_i 
    = \sum_{k=1}^{K} \mat{P}_{k}^T \Qxk \bm{v}_{i,k}
    = \sum_{k=1}^{K} \mat{P}_{k}^T \Hk \bm{v}_{i,k}'\\
    &= \left[ \sum_{k=1}^{K} \mat{P}_k^T \Hk \mat{P}_k \right] \vgd_i' 
    = \tM \vgd_i'.
\end{align*}
In the steps leading to this result, we used the fact that the generic SBP operator $\Dxk$ exactly differentiates polynomials of total degree $p$; specifically, we used the equivalent statement $\Qxk \bm{v}_{i,k} = \Hk \bm{v}_{i,k}'$.  We also used the exactness of the prolongation operators when applied to the polynomials $\vgd_i$ and $\vgd_i'$; see Remark~\ref{rmk:prolong_direct}.

Next, we need to show that the mass, or norm, matrix is symmetric positive definite.  Symmetry is obvious from the definition. The mass matrix is also positive definite; if $\ugd \in \mathbb{R}^K$ is an arbitrary vector, then
\begin{equation*}
    \ugd^T \tM \ugd = \ugd^T \left[ \sum_{k=1}^K \mat{P}_k^T \Hk \mat{P}_{k} \right] \ugd 
    = \sum_{k=1}^K \bm{u}_{k}^T \Hk \bm{u}_{k}
    = \sum_{k=1}^K \sum_{q=1}^{n_q} (\Hk)_{qq} (\bm{u}_k)_{q}^2.
\end{equation*}
This sum is clearly non-negative, since the diagonal entries $(\Hk)_{qq}$ are strictly positive for diagonal-norm SBP operators.  Furthermore, the sum is strictly positive if $\ugd$ is nonzero.  To see this, suppose otherwise; that is, suppose $\ugd^T \tM \ugd = 0$ for some $\ugd \neq \bm{0}$.  Then we must have $\sum_{k=1}^K \uk^T \Hk \uk = 0$, which is only possible if $\uk = \bm{0},\,\forall\,k$; in other words
\begin{equation*}
    \begin{bmatrix} \mat{P}_1 \\ \mat{P}_2 \\ \vdots \\ \mat{P}_K \end{bmatrix} \ugd = \begin{bmatrix} \bm{u}_1 \\ \bm{u}_2 \\ \vdots \\ \bm{u}_K \end{bmatrix} = \begin{bmatrix} \bm{0} \\ \bm{0} \\ \vdots \\ \bm{0} \end{bmatrix}.
\end{equation*}
This contradicts the assumption that the matrix on the left has full rank of $K$.  Thus, we have shown that $\tM$ is symmetric positive definite.

We also need to show that the symmetric part of $\tQx$ is equal to $\frac{1}{2}\tEx$.  This follows easily from Lemma~\ref{lem:DGD_SBP_matrices}, since
\begin{equation*}
    \tQx + \tQx^T = \sum_{k=1}^K \mat{P}_k^T \underbrace{(\Qxk + \Qxk^T)}_{\Exk} \mat{P}_k
    = \tEx.
\end{equation*}

Finally, we need to show that $\tEx$ satisfies \eqref{eq:SBP_Ex_accuracy}.  As before, let $\vgd_i$ and $\vgd_j$ denote the basis functions $\fnc{V}_i$ and $\fnc{V}_j \in \{ \fnc{V}_r \}_{r=1}^{n_p}$, respectively, evaluated at the centroids of the elements.  Then we have
\begin{align*}
    \vgd_i^T \tEx \vgd_j 
    = \vgd_i^T \left[\sum_{k=1}^K \mat{P}_k^T \Exk \mat{P}_k \right] \vgd_j 
    &= \sum_{k=1}^K \bm{v}_{i,k}^T \Exk \bm{v}_{j,k} \\
    &= \sum_{k=1}^K \int_{\Gamma_k} \fnc{V}_i \fnc{V}_j n_x \, d\Gamma 
    = \int_{\Gamma} \fnc{V}_i \fnc{V}_j n_x \, d\Gamma.
\end{align*}
We arrived at the final line by using the fact that diagonal-norm SBP operators satisfy \eqref{eq:SBP_Ex_accuracy}, and the fact that the surface integrals over interior faces cancel for polynomials of degree $p$.
\qed
\end{proof}

\ignore{
While the DGD matrix $\tDx$ satisfies the multidimensional SBP operator definition from~\cite{hicken:mdimsbp2016}, its corresponding boundary operator $\tEx$ is not typical of most SBP operators described in the literature.  Specifically, while $\tEx$ satisfies \eqref{eq:SBP_Ex_accuracy}, it does not mimic boundary integration in a discrete sense.  For example, for the arbitrary states $\ugd$ and $\vgd \in \mathbb{R}^{K}$ (\ie not polynomials of degree $p$), we have
\begin{equation*}
    \vgd^T \tEx \ugd = \sum_{k=1}^{K} \int_{\Gamma_k} \fnc{V}_h \fnc{U}_h n_x \, d\Gamma 
    = \sum_{\gamma \in \Gamma_I} \int_{\gamma} \llbracket \fnc{V}_h \fnc{U}_h \rrbracket n_x \,d\Gamma + \sum_{\gamma \in \Gamma_B} \int_{\gamma} \fnc{V}_h \fnc{U}_h n_x \, d\Gamma,
\end{equation*}
where $\llbracket \fnc{V}_h \fnc{U}_h  \rrbracket$ is the jump in the product $\fnc{V}_h \fnc{U}_h$ across the interface $\gamma$, $\Gamma_I$ is the set of element interfaces, and $\Gamma_B$ is the set of element boundaries.  For unstructured-mesh DGD schemes, the product $\fnc{V}_h \fnc{U}_h$ is not continuous across element boundaries, so it follows that the first sum on the right-hand side does not vanish.  Thus, in contrast with conventional SBP schemes, $\vgd^T \tEx \ugd$ has contributions from all elements, in general, and not merely elements adjacent to the domain boundary $\Gamma$.
}

\section{Entropy-stable discontinuous Galerkin difference discretizations}
\label{sec:entropytheory}

Having established the connection between Galerkin difference and SBP operators, we can exploit existing entropy-stable SBP theory to construct entropy-stable DGD schemes.  To that end, this section begins with a brief review of conservative hyperbolic systems that admit a strongly convex entropy function.  We then define semi-discrete entropy conservation and stability in the context of generic diagonal-norm SBP methods.  Subsequently, we show how these SBP methods can be used to construct entropy-conservative/stable DGD spatial discretizations.  We conclude by describing our entropy-stable temporal discretization.

\subsection{Hyperbolic conservation laws with convex entropy}

Consider a generic, hyperbolic conservation law in two dimensions, given by the PDE 
\begin{equation}\label{eq:cons_law}
\frac{\partial \bfnc{U}}{\partial t} + \frac{\partial \bfnc{F}_x}{\partial x} + \frac{\partial \bfnc{F}_y}{\partial y} = \bm{0}, \quad \forall \bm{x} \in \Omega, \quad \forall t \in [0,T],
\end{equation}
where $\fnc{U}(\bm{x},t) \in \mathbb{R}^{s}$ is the vector of $s$ conservative variables, and $\bfnc{F}_{x} : \mathbb{R}^s \rightarrow \mathbb{R}^s$ and $\bfnc{F}_{y} : \mathbb{R}^s \rightarrow \mathbb{R}^s$ are smooth flux functions in the $x$ and $y$ coordinate directions, respectively.  

We narrow our focus to conservation laws that have an associated convex entropy function, $\fnc{S} : \mathbb{R}^s \rightarrow \mathbb{R}$, that satisfies
\begin{equation}\label{eq:symmetrize}
    \frac{\partial^2 \fnc{S}}{\partial \bfnc{U}^2} \frac{\partial \bfnc{F}_{*}}{\partial \bfnc{U}} =
    \left[ \frac{\partial^2 \fnc{S}}{\partial \bfnc{U}^2} \frac{\partial \bfnc{F}_{*}}{\partial \bfnc{U}} \right]^T,
    \qquad \bfnc{F}_{*} \in \{ \bfnc{F}_x, \bfnc{F}_y \},
\end{equation}
where $[\partial \bfnc{F}_{*}/\partial \bfnc{U}]_{ij} = \partial \fnc{F}_{*,i}/\partial \fnc{U}_j$ is a flux Jacobian, and $[\partial^2 \fnc{S}/ \partial \bfnc{U}]_{ij} = \partial \fnc{S} / \partial \fnc{U}_i \partial \fnc{U}_j$ is the positive definite Hessian of the entropy.  The relations~\eqref{eq:symmetrize} imply~\cite{tadmor:2003} the existence of entropy fluxes $\fnc{G}_x : \mathbb{R}^s \rightarrow \mathbb{R}$ and $\fnc{G}_y : \mathbb{R}^s \rightarrow \mathbb{R}$, whose gradients obey\footnote{We follow the convention that gradients are row vectors}

\begin{equation}\label{eq:entropy_flux}
    \frac{\partial \fnc{G}_{*}}{\partial \bfnc{U}} = \bfnc{W}^T \frac{\partial \bfnc{F}_{*}}{\partial \bfnc{U}},
\end{equation}
where we have introduced the entropy variables, $\bfnc{W} \equiv [\partial \fnc{S} / \partial \bfnc{U}]^T$, for convenience.  For smooth states,
\eqref{eq:cons_law} and \eqref{eq:entropy_flux} imply that the entropy is conserved:
\begin{equation}\label{eq:entropy_cons}
    \int_{\Omega} \bfnc{W}^T \left[\frac{\partial \bm{\mathcal{U}}}{\partial t} + \frac{\partial \bfnc{F}_x}{\partial x} + \frac{\partial \bfnc{F}_y}{\partial y} \right] \,d\Omega
    = \int_{\Omega} \left[ \frac{\partial \fnc{S}}{\partial t} +
    \frac{\partial \fnc{G}_x}{\partial x} + \frac{\partial \fnc{G}_y}{\partial y} \right] \, d\Omega = 0.
\end{equation}
More generally, the entropy for physically relevant weak solutions to \eqref{eq:cons_law} satisfies the following inequality (refer to the review~\cite{tadmor:2003}, and the references therein):
\begin{equation}\label{eq:entropy_stable}
\int_{\Omega} \left[ \frac{\partial \fnc{S}}{\partial t} +
    \frac{\partial \fnc{G}_x}{\partial x} + \frac{\partial \fnc{G}_y}{\partial y} \right] \, d\Omega \leq 0.
\end{equation}

\subsection{Entropy-conservative and entropy-stable SBP discretizations}

We are interested in discretizations that mimic \eqref{eq:entropy_cons} and \eqref{eq:entropy_stable}, because these properties imply an $L^2$ bound on the state~\cite{Dafermos2010hyperbolic}, and, consequently, they impart a form of nonlinear stability.  What we mean by mimic will be made precise for diagonal-norm SBP discretizations in Definition~\ref{def:entropy_dis} below.

Consider the following SBP discretization of \eqref{eq:cons_law} defined on each element $k$ in the tesselation $T_h$ of $\Omega$.
\begin{equation}\label{eq:SBP_entstable}
    \frac{d \uk}{dt} + \bm{r}_k(\bm{u}) = \bm{0},
    \qquad\forall k = 1,2,\ldots,K,
\end{equation}
where $\uk \in \mathbb{R}^{sn_q}$ holds the conservative variables at the quadrature nodes $X_k = \{ \bm{x}_q \}_{q=1}^{n_q}$ on element $k$, and $\bm{u} = [\bm{u}_1^T, \bm{u}_2^T, \ldots, \bm{u}_K^T ]^T \in \mathbb{R}^{s n_q K}$ is the compound vector of conservative variables over all $n_q K$ nodes.  We assume the SBP operators used in the spatial discretization $\bm{r}_k$ are associated with a diagonal norm matrix $\Hk \in \mathbb{R}^{n_q \times n_q}$.  Note that $\bm{r}_k(\bm{u})$ depends on the neighbours of element $k$ through interface fluxes, which is why $\bm{r}_k$ is not a function of $\uk$ alone.

In order to define discrete entropy conservation and stability for \eqref{eq:SBP_entstable}, we need to be able to relate the discrete conservative and entropy variables to one another.  To this end, we will assume that the unknowns in $\uk$ are ordered such that the conservative variables at an arbitrary node $\bm{x}_q \in X_k$ are sequential:
\begin{equation*}
    \bm{u}_q = [\uk]_{s(q-1)+1:sq}, \qquad \forall\, q = 1,2,\ldots,n_q.
\end{equation*}
Then, it follows that the discrete entropy variables at the nodes are given by
\begin{equation*}
    \bm{w}_q = \bfnc{W}(\bm{u}_q), \qquad \forall\, q = 1,2,\ldots,n_q.
\end{equation*}
We will use $\wk = [\bm{w}_1^T, \bm{w}_2^T,\ldots, \bm{w}_{n_q}^T]^T$ to denote the set of all entropy variables on element $k$.

We can now define what we mean by entropy conservative and entropy stable SBP discretizations.
\begin{definition}[Entropy-conservative/stable SBP discretizations]\label{def:entropy_dis}
 The spatial discretization \eqref{eq:SBP_entstable} is entropy-conservative with respect to the diagonal, positive definite matrices $\Hk \in \mathbb{R}^{n_q \times n_q}$, $k=1,2,\ldots,K$, if
\begin{equation}\label{eq:ent_cons_cond}
    \sum_{i=k}^{K} \wk^T \Hk \bm{r}_k(\bm{u}) = 0.
\end{equation}
Similarly, \eqref{eq:SBP_entstable} is entropy-stable with respect to the $\Hk$ if
\begin{equation}\label{eq:ent_stable_cond}
\sum_{i=1}^{n} \wk^T \Hk \bm{r}_k(\bm{u}) \geq 0.
\end{equation}
\end{definition}

\begin{remark}
Definition~\ref{def:entropy_dis} is general enough to accommodate boundary conditions, although constructing entropy-conservative and entropy-stable boundary conditions is nontrivial, in general.
\end{remark}

From Definition~\ref{def:entropy_dis}, we can infer that an entropy-conservative SBP semi-discretization satisfies
\begin{equation}\label{eq:entropy_cons_dis}
    \sum_{i=k}^{K} \wk^T \Hk \left[ \frac{d \uk}{dt} + \bm{r}_k(\bm{u}) \right] 
    = \sum_{k=1}^{K} \sum_{q=1}^{n_q} (\Hk)_{qq} \bm{w}_{q}^T  \frac{d \bm{u}_q}{dt} = \sum_{k=1}^{K} \sum_{q=1}^{n_q} (\Hk)_{qq} \frac{d s_q}{dt} = 0,
\end{equation}
where $s_q= \fnc{S}(\bm{u}_q)$ is the entropy based on $\bm{u}_q$.  Similarly, an entropy-stable SBP discretization implies the following condition:
\begin{equation}\label{eq:entropy_stable_dis}
    \sum_{k=1}^{K} \sum_{q=1}^{n_q} (\Hk)_{qq} \frac{d s_q}{dt} \leq 0.
\end{equation}

On periodic domains, \eqref{eq:entropy_cons_dis} and \eqref{eq:entropy_stable_dis} are discrete analogs of \eqref{eq:entropy_cons} and \eqref{eq:entropy_stable}, respectively, since the divergence of $(\fnc{G}_x,\fnc{G}_y)$ will vanish on a periodic domain.  Furthermore, since the SBP nodes $X_k$ and corresponding diagonal entries in $\Hk$ define a quadrature rule that is degree $2p-1$ exact (at least), then
\begin{equation*}
    \sum_{k=1}^{K} \sum_{q=1}^{n_q} (\Hk)_{qq} \frac{d s_q}{dt} = \int_{\Omega} \frac{\partial \fnc{S}}{\partial t} \, d\Omega + \text{O}(h^{2p}).
\end{equation*}

\subsection{Construction of entropy-conservative/stable DGD discretizations from diagonal-norm SBP discretizations}

Entropy-conservative and entropy-stable DGD discretizations can be constructed from SBP discretizations following an approach similar to Theorem~\ref{thm:DGD_SBP_equiv} for constant-coefficient linear advection.  That is, we can define the solution at the SBP nodes by prolonging the DGD solution.  However, unlike the linear advection case, the choice of \emph{which} variable to prolong is important in the context of nonlinear conservation laws.  We must prolong the entropy variables, not the conservative variables, to ensure the DGD scheme is entropy conservative/stable.

Let $\wgd \in \mathbb{R}^{sK}$ denote the entropy variables associated with the DGD degrees of freedom at the center of each element.  As with the SBP discretization, suppose the variables are ordered such that entropy variables at the $k$th degree of freedom are sequential:
\begin{equation*}
    \wgd_k = [\wgd]_{s(k-1)+1:sk},
    \qquad\forall\, k=1,2,\ldots,K.
\end{equation*}
Note that the DGD vector $\wgd_k$ should not be confused with the SBP vector $\wk$: the former has $s$ entries while the latter has $sn_q$ entries.

With the above ordering, we can prolong the entropy variables to the SBP nodes of element $k$ using the (scalar) prolongation operator $\Pk$, and the Kronecker product $\otimes$:
\begin{equation*}
    \wk = (\Pk \otimes \mat{I}_s) \wgd = \barPk \wgd ,
    \qquad \forall\, k=1,2,\ldots,K,
\end{equation*}
where $\mat{I}_s$ is the $s\times s$ identity matrix, and $\barPk \equiv \Pk \otimes \mat{I}_s$.  Thus, for the DGD discretization, the entropy variables at the SBP nodes are given by
\begin{equation*}
    \bm{w}_q = \bm{w}_q(\wgd) = \sum_{\nu=1}^{K} (\Pk)_{q\nu} \wgd_{\nu},
    \qquad \forall\, \bm{x}_q \in X_k,\quad k=1,2,\ldots,K.
\end{equation*}
Finally, the conservative variables at $\bm{x}_q$ are obtained by transforming the entropy variables $\bm{w}_q$, \viz 
\begin{equation*}
\bm{u}_q = \bfnc{U}(\bm{w}_q),\qquad \forall\, q=1,2,\ldots,n_q,
\end{equation*}
and we will use the compound vector $\uk(\wgd) = [\bm{u}_1^T, \bm{u}_2^T, \ldots, \bm{u}_{n_q}^T]^T$  to denote the conservative variables at the SBP nodes of element $k$ computed from the prolonged entropy variables.

We construct the DGD discretization by substituting $\uk(\wgd)$ into~\eqref{eq:SBP_entstable}, left multiplying each equation by $\barPk^T \barMk$, and then summing the result over all elements:
\begin{equation}\label{eq:DGD_entstable}
    \sum_{k=1}^{K} \barPk^T \barMk \left[ \frac{d}{dt}\uk(\wgd) + \bm{r}_k(\wgd)\right] = \bm{0},
\end{equation}
where $\bm{r}_k(\wgd)$ is the SBP spatial discretization on element $k$ evaluated using $\uk(\wgd)$, $k=1,2,\ldots,K$.  The discretization~\eqref{eq:DGD_entstable} is a system of $sK$ nonlinear ODEs in $sK$ unknowns (\ie the entries in $\wgd$).

\begin{thrm} \label{thm:DGD_entstable}
If the SBP discretization~\eqref{eq:SBP_entstable} is entropy-stable, then the DGD discretization~\eqref{eq:DGD_entstable} satisfies
\begin{equation*}
 \sum_{k=1}^{K} \sum_{q=1}^{n_q} (\Hk)_{qq} \frac{d s_q}{dt} \leq 0,
\end{equation*}
where $s_q = \fnc{S}(\bm{w}_q(\wgd))$ denotes the entropy at SBP node $q$ of element $k$ based on the prolonged DGD entropy variables.  The inequality becomes an equality if the SBP spatial discretization is entropy conservative.
\end{thrm}

\begin{proof}
Suppose \eqref{eq:SBP_entstable} is entropy stable.  Then, if we left multiply the DGD discretization \eqref{eq:DGD_entstable} by $\wgd^T$ and recall that $\wk = \barPk \wgd$, we find
\begin{multline}\label{eq:DGD_entstable_eq1}
\sum_{k=1}^{K} \wgd^T \barPk^T \barMk \left[ \frac{d}{dt}\uk(\wgd) + \bm{r}_k(\wgd)\right] \\
= \sum_{k=1}^{K} \wk^T \barMk \frac{d}{dt} \uk(\wgd) + 
\sum_{k=1}^{K} \wk^T \barMk \bm{r}_k(\wgd) = 0.
\end{multline}
Now, the spatial term $\bm{r}_k(\wgd)$ is evaluated using the conservative variables $\uk(\wgd)$, which are themselves evaluated using the prolonged entropy variables $\wk = \barPk \wgd$.  Consequently, since \eqref{eq:SBP_entstable} is entropy stable by assumption, we have $\sum_{k=1}^K \wk^T \barMk \bm{r}_k(\wgd) \geq 0$, so  \eqref{eq:DGD_entstable_eq1} implies
\begin{equation*}
\sum_{k=1}^{K} \sum_{q=1}^{n_q} (\Hk)_{qq} \bm{w}_q^T \frac{d \bm{u}_q}{dt} 
= \sum_{k=1}^{K} \sum_{q=1}^{n_q} (\Hk)_{qq} \frac{d s_q}{dt}
\leq 0,
\end{equation*}
as desired.  The proof for the corresponding entropy-conservative spatial discretizations is similar.
\end{proof}

\begin{remark}
Note that the proof fails if we prolong the conservative variables, because the prolongation operation does not commute with the nonlinear entropy-variable transformation:
\begin{equation*}
    \bm{w}_q = \sum_{\nu=1}^K (\Pk)_{q\nu} \wgd_{\nu} \neq 
    \bfnc{W}\left({\textstyle \sum_{\nu=1}^K (\Pk)_{q\nu} \ugd_{\nu}} \right),
\end{equation*}
where $\ugd_\nu$ denotes the conservative variables at DGD degree of freedom $\nu$.
\end{remark}

\begin{remark}
While the DGD discretization must prolong the entropy variables to the SBP nodes to ensure entropy conservation or stability, this does not preclude using the conservative variables as the unknown state.  If the conservative variables are adopted as the state, it merely introduces an additional step, namely, conversion to the entropy variables, before prolongation.  Such an approach is used in \cite{Chan2018discretely} in the context of discontinuous Galerkin finite-element methods.
\end{remark}

\subsection{Entropy-stable temporal discretization}\label{implicit_rrk}

So far we have focused exclusively on semi-discrete entropy-conservation and stability, but it is desirable to mimic \eqref{eq:entropy_cons} and \eqref{eq:entropy_stable} in a fully discrete sense to ensure entropy remains bounded.  Therefore, for our temporal discretization, we use an entropy-stable relaxation Runge-Kutta (RRK) as described by Ranocha \etal~\cite{Ranocha2020relaxation}.  This section briefly summarizes our RRK implementation.

We adopt the implicit midpoint method as our baseline temporal discretization, which is a second-order accurate Gauss-Legendre method.  Thus, to move from time $t^n$ to $t^{n+1}$ we first solve the midpoint discretization of \eqref{eq:DGD_entstable}:
\begin{equation}\label{eq:DGD_midpoint}
    \sum_{k=1}^{K} \barPk^T \barMk \left[ \uk(\wgd^*) - \uk(\wgd^n) + \Delta t\, \bm{r}_k(\wgd^{n+\frac{1}{2}}) \right] = \bm{0},
\end{equation}
where $\Delta t = t^{n+1} - t^{n}$ and
\begin{equation*}
    \wgd^{n+\frac{1}{2}} \equiv \frac{1}{2}\wgd^* + \frac{1}{2}\wgd^n.
\end{equation*}
In our implementation, we solve~\eqref{eq:DGD_midpoint} for $\wgd^{*}$ using Newton's method and sparse direct solves for the linear subproblems.  The tolerances used for Newton's method will be discussed in Section~\ref{numericalexp}.

In the conventional midpoint method, we would set $\wgd^{n+1} = \wgd^{*}$.  However, even though the midpoint method is a symplectic integrator,  $\wgd^{*}$ does not necessarily ensure entropy conservation or stability.  This is because the fully discrete version of Theorem~\ref{thm:DGD_entstable} does not hold automatically, since
\begin{equation*}
    \sum_{k=1}^K (\wgd^{n+\frac{1}{2}})^T \barPk^T \barMk \left[  \uk(\wgd^*) - \uk(\wgd^n) \right] 
    \neq \sum_{k=1}^K \sum_{q=1}^{n_q} (\Hk)_{qq} \left[s_q^* - s_q^n \right],
\end{equation*}
where $s_q^{*}$ and $s_q^n$ are the nodal entropies based on $\wgd^{*}$ and $\wgd^{n}$, respectively.

Rather than setting $\wgd^{n+1} = \wgd^{*}$, the RRK version of the midpoint method sets
\begin{equation}\label{eq:RRK_update}
\wgd^{n+1}_{\beta^n} = \wgd^{n} + \beta^n(\wgd^{*} - \wgd^n),
\end{equation}
where $\beta^n \in \mathbb{R}$ is chosen such that entropy conservation/stability is respected discretely.  To be more precise, let
\begin{equation}\label{eq:total_entropy}
    S(\wgd) \equiv \sum_{k=1}^K \sum_{q=1}^{n_q} (\Hk)_{qq} s_q(\wgd) 
\end{equation}
denote the global (integrated) entropy based on the DGD entropy variables $\wgd$, and let 
\begin{equation*}
    R(\wgd) \equiv \sum_{k=1}^K \wgd^T \barPk^T \barMk \bm{r}_k(\wgd)
\end{equation*}
denote the change in entropy due to the spatial discretization.  Then we seek $\beta^n$ such that 
\begin{equation}\label{eq:RRK_residual}
    S(\wgd^{n+1}_{\beta^n}) - S(\wgd^n) + \beta^n \Delta t R(\wgd^{n+\frac{1}{2}}) = 0,
\end{equation}
where $\wgd^{n+1}_{\beta^n}$ is defined by \eqref{eq:RRK_update}. Note that \eqref{eq:RRK_residual} is a scalar equation in the unknown $\beta^n$, which we solve to a tolerance of $10^{-13}$ using the secant method.  Once $\beta^n$ has been determined, $\wgd^{n+1}$ is computed using \eqref{eq:RRK_update}.

\begin{remark}\label{remark:explicit}
The entropy-stable DGD discretization under consideration does not appear to be well suited to explicit RRK schemes.  For example, the forward Euler scheme would require solving
\begin{equation*}
    \sum_{k=1}^{K} \barPk^T \barMk \left[ \uk(\wgd^*) - \uk(\wgd^n) + \Delta t\, \bm{r}_k(\wgd^n) \right] = \bm{0}
\end{equation*}
for $\wgd^{*}$.  Since $\uk(\wgd^*)$ is a nonlinear function of $\wgd^{*}$, the solution cannot be determined without using Newton's method.
\end{remark}

\section{Numerical experiments} \label{numericalexp}

This section presents numerical experiments to verify the accuracy and stability of a particular the DGD implementation.  We also use these experiments to investigate the spectra of entropy-stable DGD discretizations.

The section begins with a description of our particular DGD implementation.  We then review the Euler equations of gas dynamics, since they are the hyperbolic conservation law targeted by the numerical experiments.  The numerical results are presented last.

\subsection{DGD implementation including stencil construction}\label{sec:results}

Our DGD and SBP discretizations are implemented in the Modular Finite Element Methods (MFEM)~\cite{mfem-library} library, and we restrict our focus to one-dimensional intervals and two-dimensional triangular grids.  The element-level SBP discretization for $\bm{r}_k(\wgd)$, which also appears in the DGD discretization~\eqref{eq:DGD_entstable}, is based on the one presented in Reference~\cite{crean:entropy2017}.  However, we use the diagonal-norm SBP simplex operators described in~\cite{Hicken2020entropy} for the present results.

The discussion of the DGD basis in section \ref{sec:basis}, while general, omitted details regarding the construction of the stencil, or patch, $N_k = \{\nu_1,\nu_2,\ldots,\nu_{n_k}\}$, \ie, the elements that influence the solution on element $k$.  To construct $N_k$ on the triangular grids considered in this work, we follow the approach in~\cite{Li2019}, which we briefly review here.

\begin{figure}[tbp]
\centering
  \subfigure[Stencil for $p=0$]{%
    \includegraphics[width=0.32\textwidth]{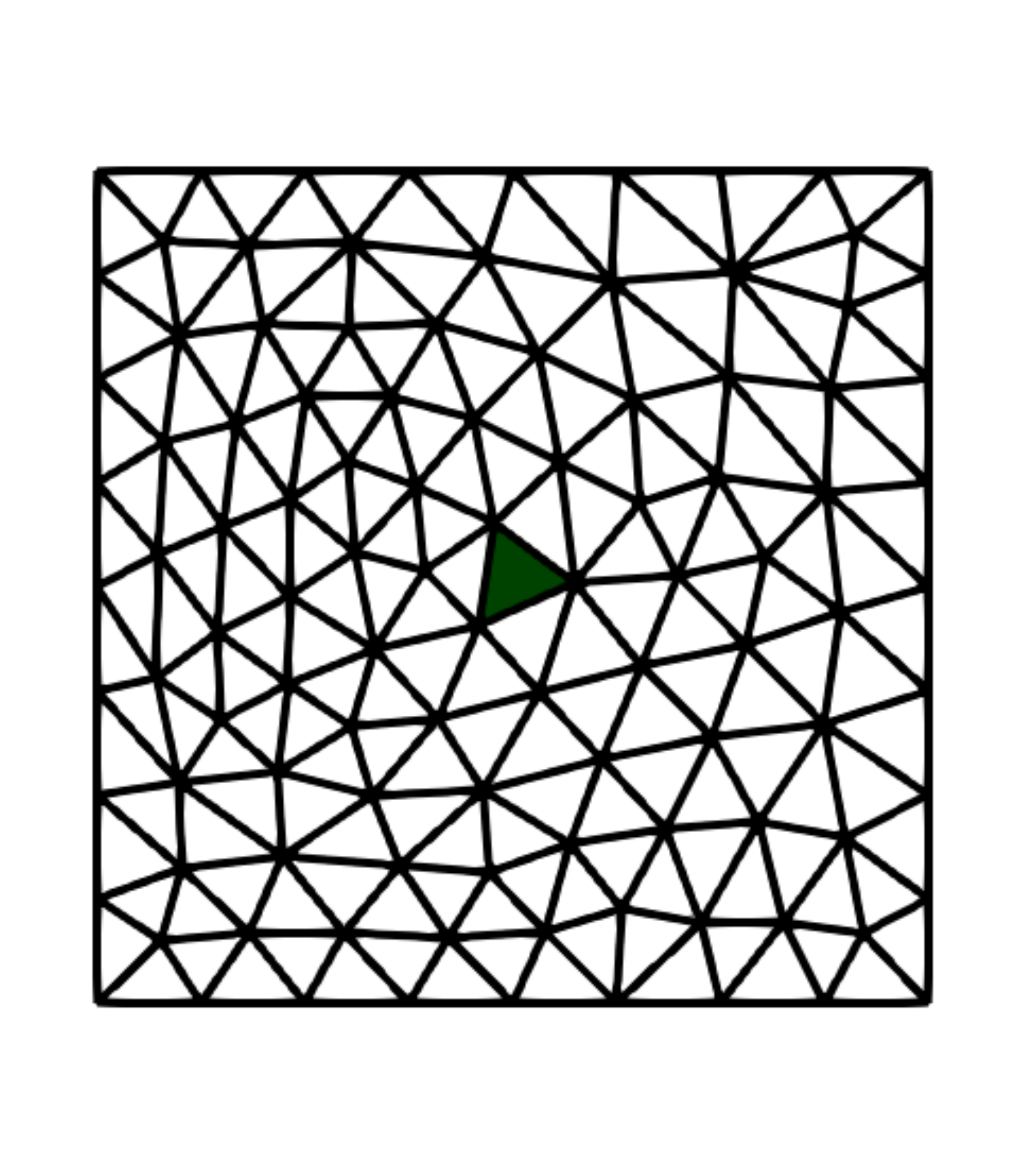}}
  \subfigure[Stencil for $p=1$]{%
    \includegraphics[width=0.32\textwidth]{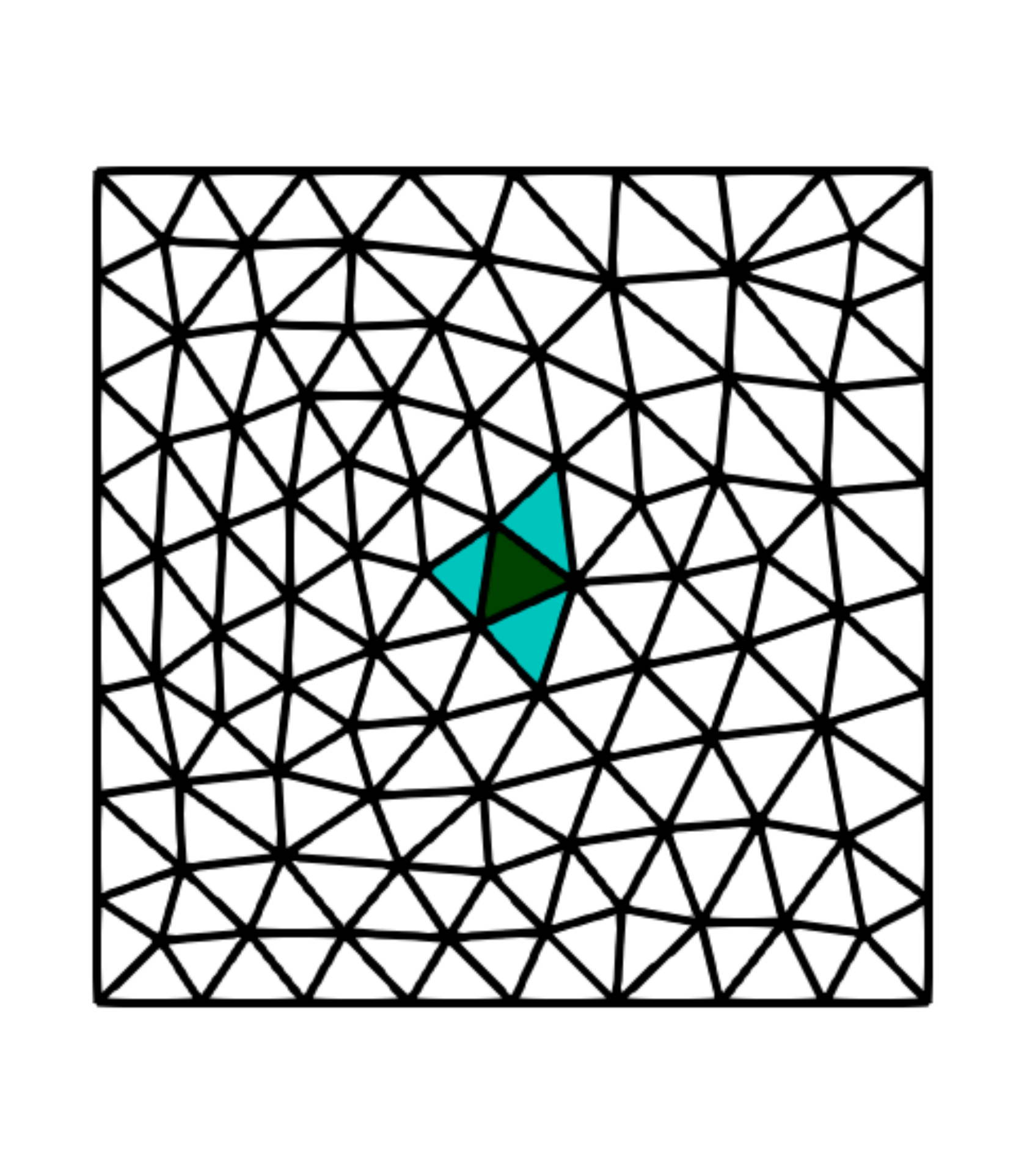}}
  \subfigure[Stencil for $p=2$]{%
    \includegraphics[width=0.32\textwidth]{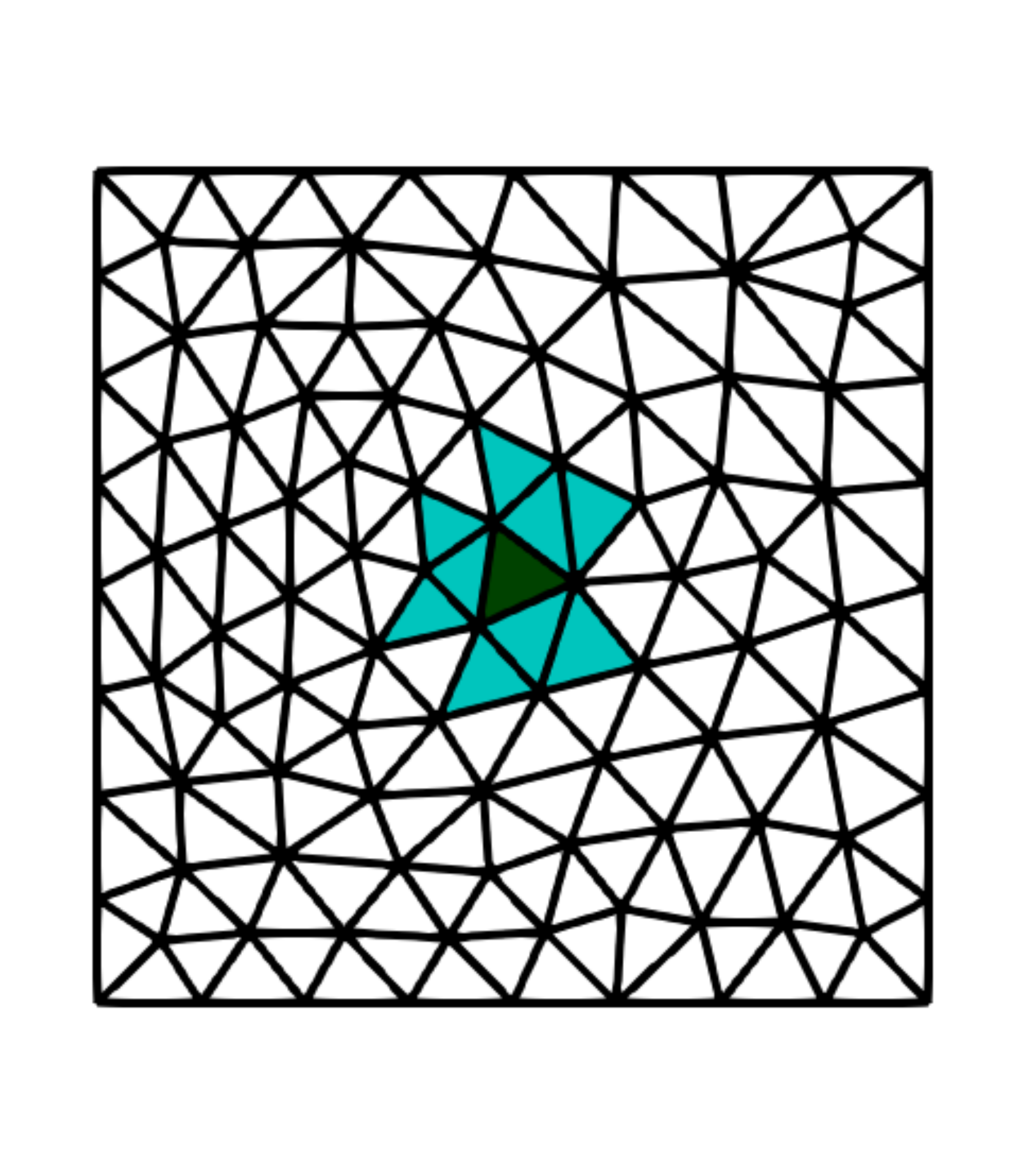}}
  \caption{Stencil construction}
\label{fig:patch}
\end{figure}

The element stencil $N_k$ is generated using an iterative process.  The set is initiated with the element $k$ itself, and then neighboring elements are added until the stencil size is at least as large as the number of basis functions for $\mathbb{P}_p(\Omega_k)$.  To be more precise, we introduce the following recursive set definition:
\begin{equation*}
\begin{split}
    N_{k}^{0} &= \{ k\}, \\
    N_{k}^{j} &= N_{k}^{j-1} \cup F(N_{k}^{j-1}),
\end{split}
\end{equation*}
where $F(N_{k}^{j-1})$ are the face-adjacent elements of the elements in $N_k^{j-1}$ in the sense of graph theory.
Then $N_k \equiv N_k^j$ for the smallest $j$ such that the number of elements in $N_k^j$ is greater than or equal to $n_p$.  Figure \ref{fig:patch} illustrates stencils that result when degree $p=0$, $p=1$, and $p=2$ basis functions are used.

\begin{remark}
Other methods can be used to define $N_k$.  For example, a stencil/patch may be defined by all elements whose centers are inside the circle of radius $r$ centered at the $k$th element's centroid.
\end{remark}

\subsection{The Euler equations}

We verify and investigate our entropy-stable DGD framework using the Euler equations.  The two-dimensional Euler equations take the form of~\eqref{eq:cons_law}, with the conservative variables given by
\begin{equation*}
 \bfnc{U} = \left[ \rho, \rho u, \rho v, e \right]^T
\end{equation*}
where $\rho$ is the density, $[\rho u, \rho v]^T$ is the momentum per unit volume, and $e$ is the total energy per unit volume.  The Euler fluxes are 
\begin{equation}\label{eq:euler_fluxes}
\bfnc{F}_x =
\begin{bmatrix}
\rho u\\
\rho u^2 + p\\
\rho u v\\
(e + p) u
\end{bmatrix},
\qquad\text{and}\qquad
\bfnc{F}_y =
\begin{bmatrix}
\rho v\\
\rho vu\\
\rho v^2 +p\\
(e + p) u
\end{bmatrix}.
\end{equation}
The pressure, $p$, is evaluated using the calorically perfect gas law as $p = (\gamma -1) \left[ e - \rho(u^2 + v^2)/2\right]$, with $\gamma = 1.4$.

Many discretizations of the Euler equations use the fluxes \eqref{eq:euler_fluxes} directly, but entropy-conservative and entropy-stable SBP discretizations rely on two-point entropy-conservative flux functions.  In this work we use the Ismail-Roe flux function~\cite{Ismail2009affordable} in our element-level discretization $\bm{r}_k$.  Thus, based on this choice, the entropy function associated with our discretization is
\begin{equation*}
 \fnc{S} = -\rho s/(\gamma -1),
\end{equation*}
where $s = \ln{(p/\rho^\gamma)}$ is the thermodynamic entropy.

\subsection{Spatial accuracy verification using the steady isentropic vortex}

We begin with a verification of spatial accuracy.  In order to isolate the spatial errors from the temporal errors, we favor a steady problem such as the isentropic vortex.
The vortex flow consists of circular streamlines and a radially varying density and pressure.  While simple, the isentropic vortex has the advantage of offering an exact solution; see, for example, \cite{hicken:dual2014}. 

\begin{figure}[tbp]
\centering
    \subfigure[Example 2D mesh for N = 10\label{fig:mesh1}]{%
    \includegraphics[width=0.46\textwidth]{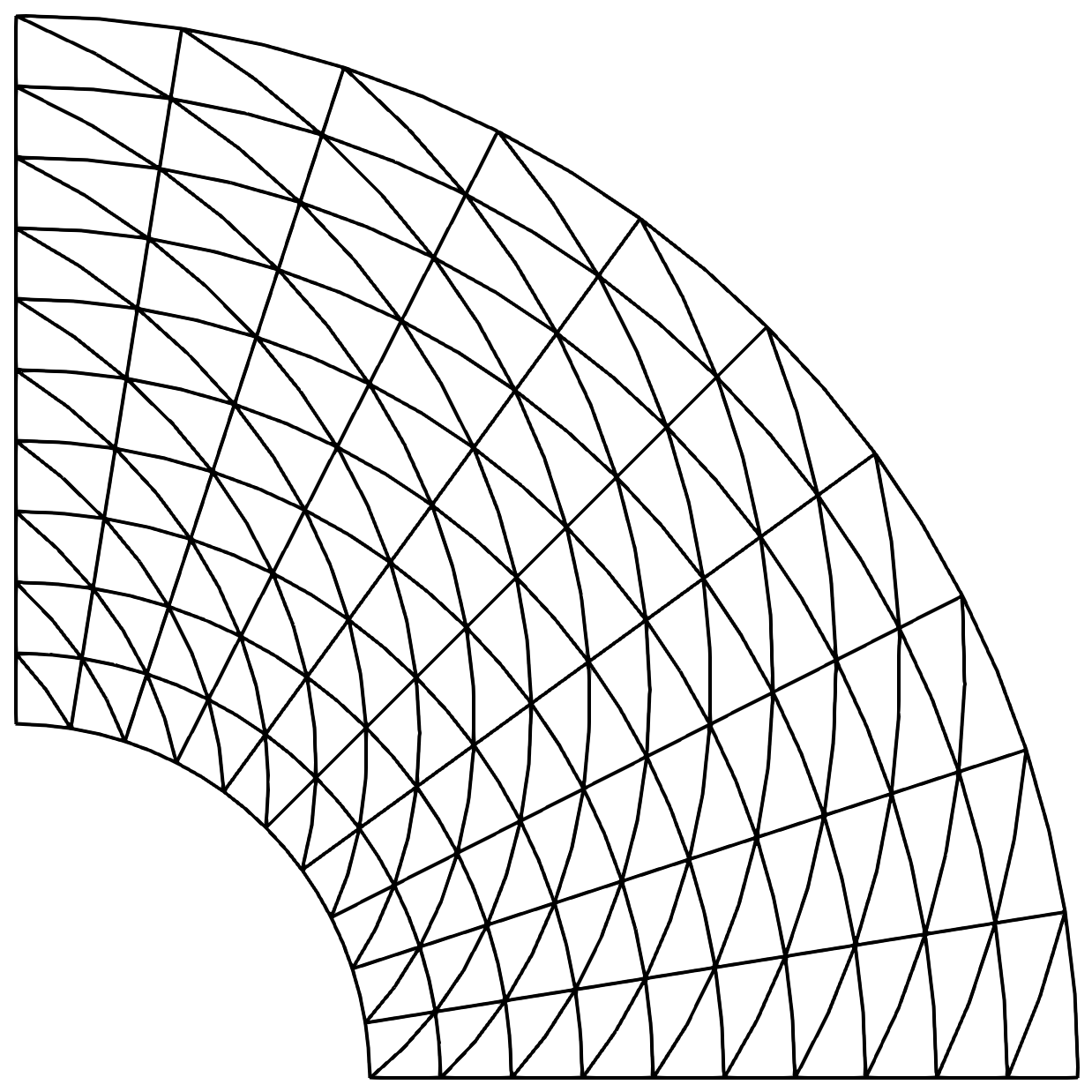}}
    \subfigure[density for $p=3$, $N= 10$\label{fig:density1}]{%
    \includegraphics[width=0.48\textwidth]{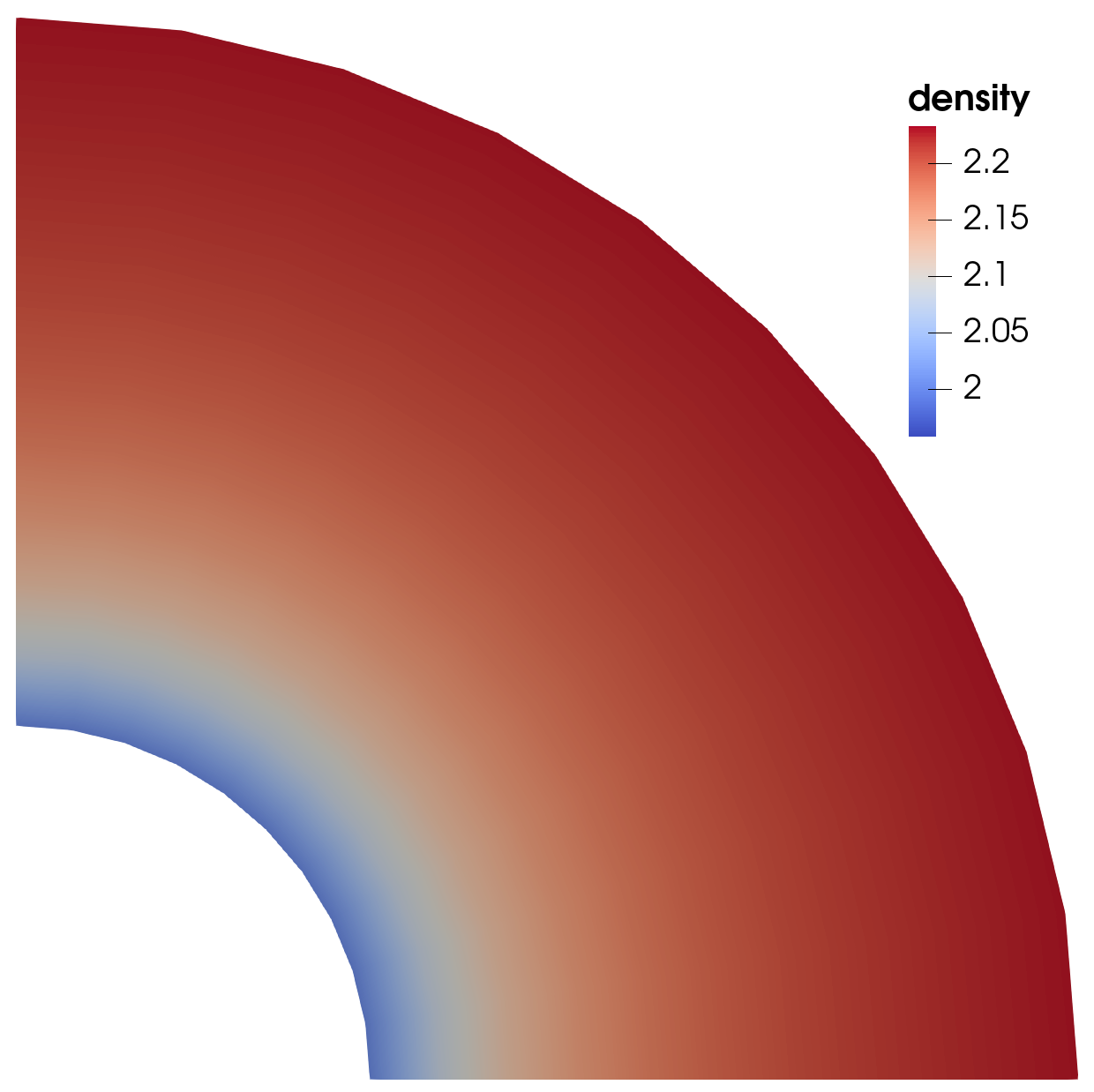}}
    \caption{Example mesh and solution for the steady-vortex problem}
\end{figure}

The isentropic vortex problem is defined on a quarter annulus domain $\Omega = \{(r, \theta), | 1 \leq r \leq 3, \; 0 \leq \theta \leq \frac{\pi}{2} \}$.  The mesh for the domain is generated in the manner described in~\cite{Hicken2020entropy}.  First, a uniform triangular mesh is defined in polar-coordinate space by splitting $N \times N$ quadrilaterals into $2N^2$ triangles.  Subsequently, the triangles are mapped to physical space using a degree $p+1$ mapping for a degree $p$ discretization. Figure~\ref{fig:mesh1} shows a sample mesh for $N = 10$, and Figure~\ref{fig:density1} illustrates the numerical density prolonged to SBP quadrature points obtained on this mesh using a $p=3$ SBP entropy-stable discretization.


As in Reference~\cite{Hicken2020entropy}, we apply a slip-wall boundary condition along the inner radius $r = 1$.  This boundary condition is imposed by evaluating the Euler flux with the state's velocity projected to be perpendicular to the surface normal.  For the remaining three sides of the domain, the exact solution is imposed weakly using the Roe numerical-flux function.  The steady versions of \eqref{eq:SBP_entstable} and \eqref{eq:DGD_entstable} are solved using Newton's method with absolute and relative tolerances both set to $10^{-14}$.

\begin{figure}[tbp]
\centering
    \subfigure[SBP discretization\label{fig:dsbp_convergence}]{%
    \includegraphics[width=0.49\textwidth]{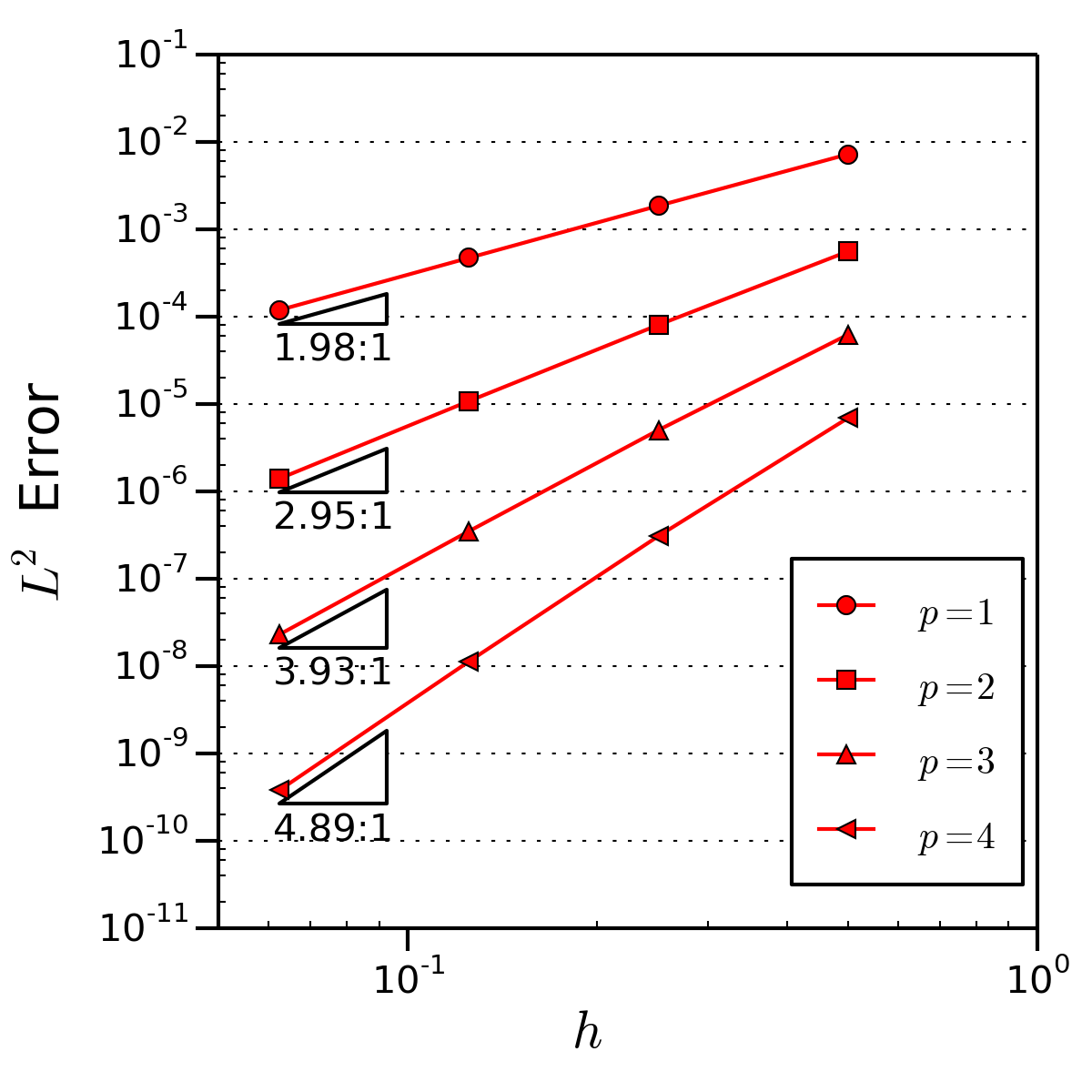}}
    \subfigure[DGD discretization\label{fig:dgd_convergence}]{%
    \includegraphics[width=0.49\textwidth]{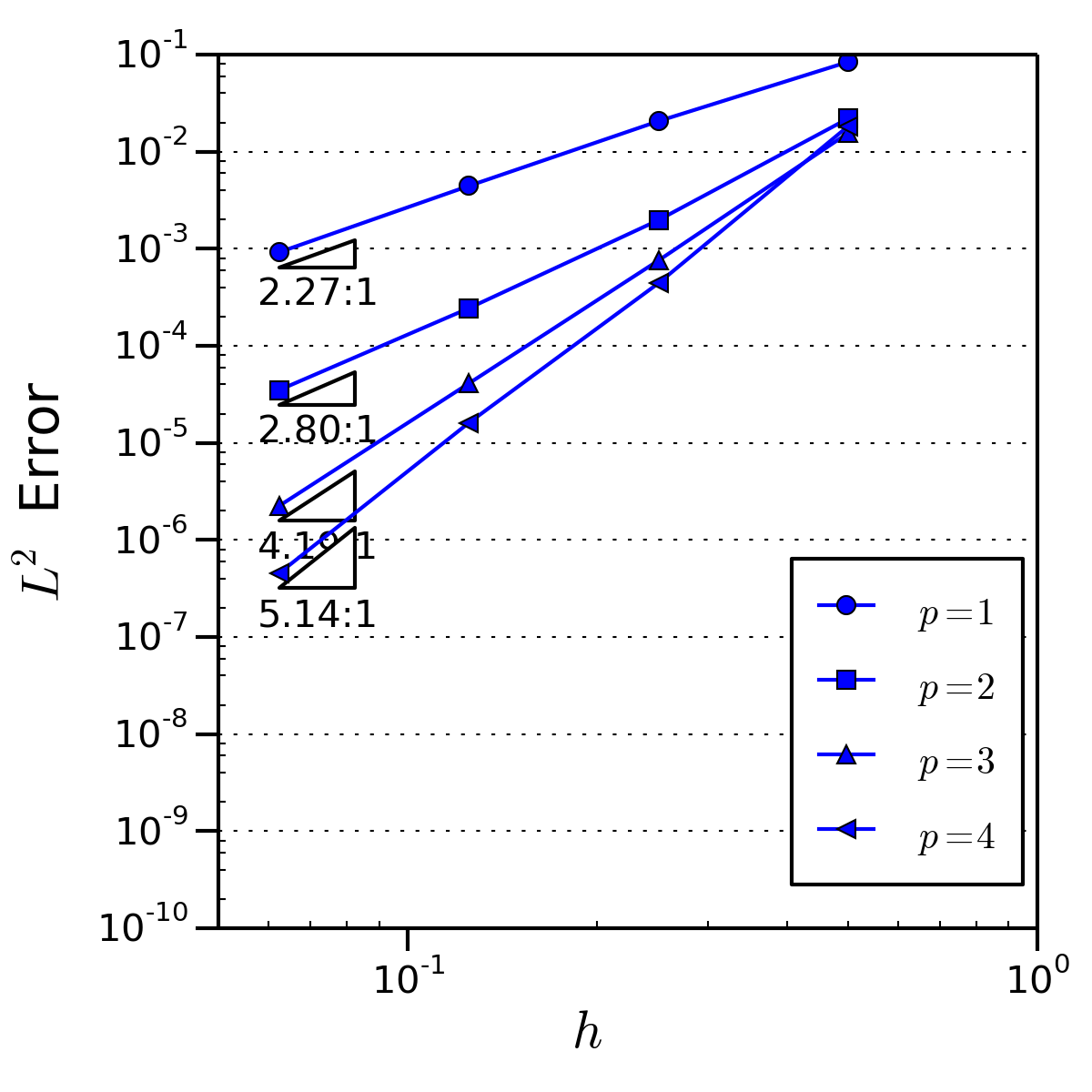}}
    \caption{Density solution $L^2$ error for SBP and DGD discretizations}
    \label{fig:solution_l2_results}
\end{figure}

Figure~\ref{fig:solution_l2_results} plots the error in the density from the SBP and DGD discretizations of degree $p=1,2,3,$ and $4$ as a function of element size, $h = 1/N$. Here, the error is an approximation to the integral $L^2$ error:
\begin{equation*}
    L^2 \; \text{Error} \approx \sqrt{\sum_{k=1}^{K} (\bm{\rho}_k - \bm{\rho}_k^{\text{exact}})^T \hat{\mat{H}}_k (\bm{\rho}_k - \bm{\rho}_k^{\text{exact}})},
\end{equation*}
where $\bm{\rho}_k$ is the density at the SBP nodes of element $k$ obtained from either a degree $p$ SBP or DGD discretization, $\bm{\rho}_k^{\text{exact}}$ is the exact density at the SBP nodes, and $\hat{\mat{H}}_k$ denotes the $2p$ exact SBP norm matrix scaled by the determinant of the mapping Jacobian.

The results in Figure~\ref{fig:dsbp_convergence} show that the errors produced by the SBP schemes approach the optimal convergence rate of $p+1$ asymptotically. The DGD rates of convergence are similar, although Figure~\ref{fig:dgd_convergence} shows the errors are somewhat super-convergent for degrees $p=1, 3, 4$ and sub-optimal for degree $p=2$.

\begin{figure}[tbp]
\centering
    \includegraphics[width=0.98\textwidth]{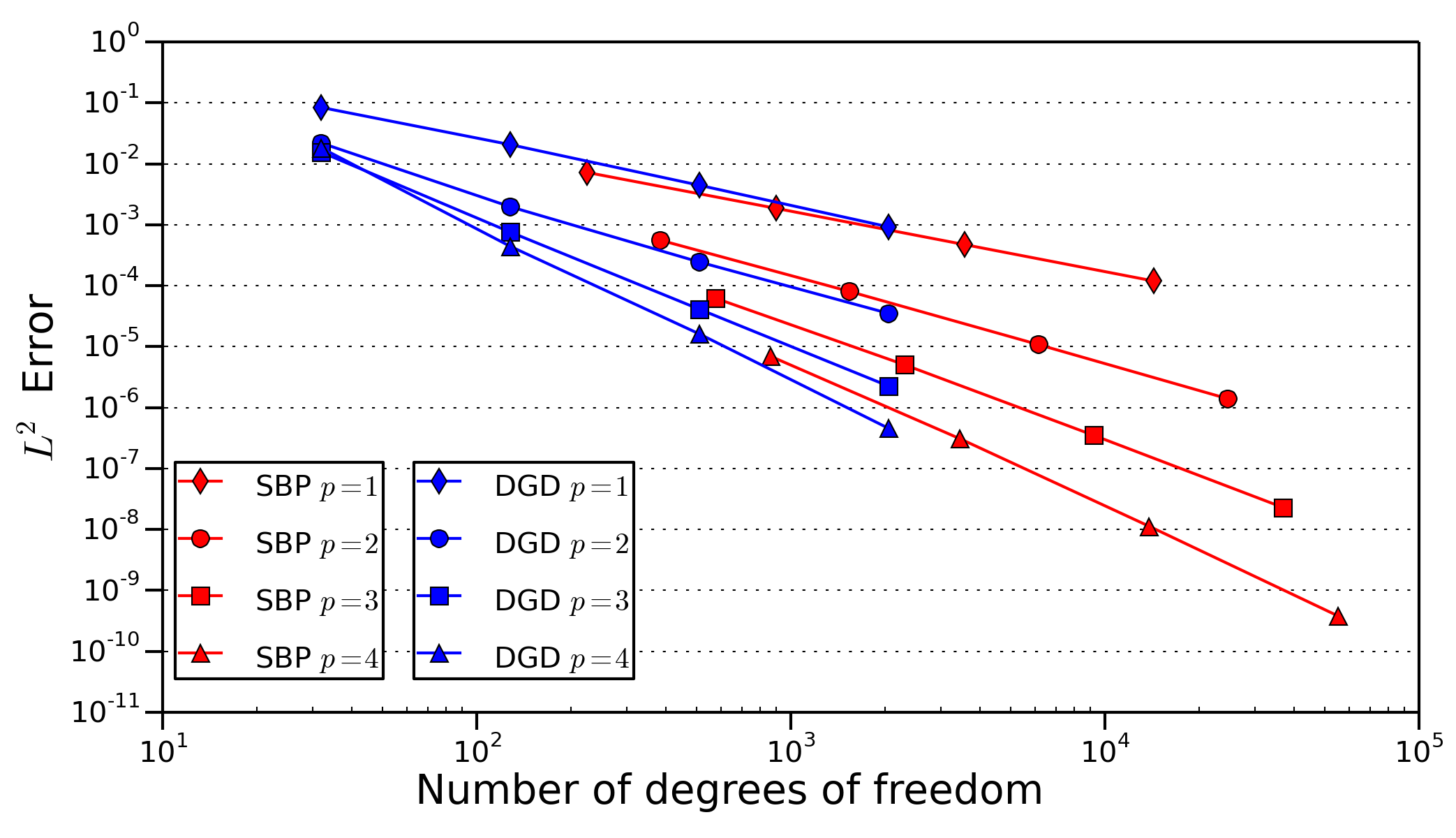}
    \caption{$L^2$ error versus degree-of-freedoms for the DGD and SBP discretizations}
    \label{fig:dofs}
\end{figure}

Comparing the SBP and DGD schemes on the basis of Figures~\ref{fig:dsbp_convergence} and \ref{fig:dgd_convergence} would be unfair, since a degree $p$ SBP scheme has significantly more degrees of freedom than a degree $p$ DGD scheme on the same mesh.  Thus, to compare these two schemes more fairly, we plot the $L^2$ error versus degrees of freedom in Figure~\ref{fig:dofs}. In this figure, the number of degrees of freedom is defined as the number of elements, $K$, for the DGD discretization and $K$ times the number of nodes per element for the SBP discretization.  Figure~\ref{fig:dofs} shows that, when measured in terms of degrees of freedom, a degree $p$ DGD discretization generally outperforms the degree $p$ SBP discretization.  The one exception is $p=1$, where the DGD error only begins to overlap the SBP error on the finest mesh.  

\begin{remark}
Error versus degrees of freedom, while better than error versus element size, remains an imperfect means of comparing the SBP and DGD schemes. In particular, the number of degrees of freedom does not reflect the potential computational saving of the DGD scheme due to its better spectral radius and conditioning; see, for example, the DGD spectra presented in the next section.
\end{remark}

\begin{figure}[tbp]
\centering
    \subfigure[SBP discretization\label{fig:dsbp_func}]{%
    \includegraphics[width=0.49\textwidth]{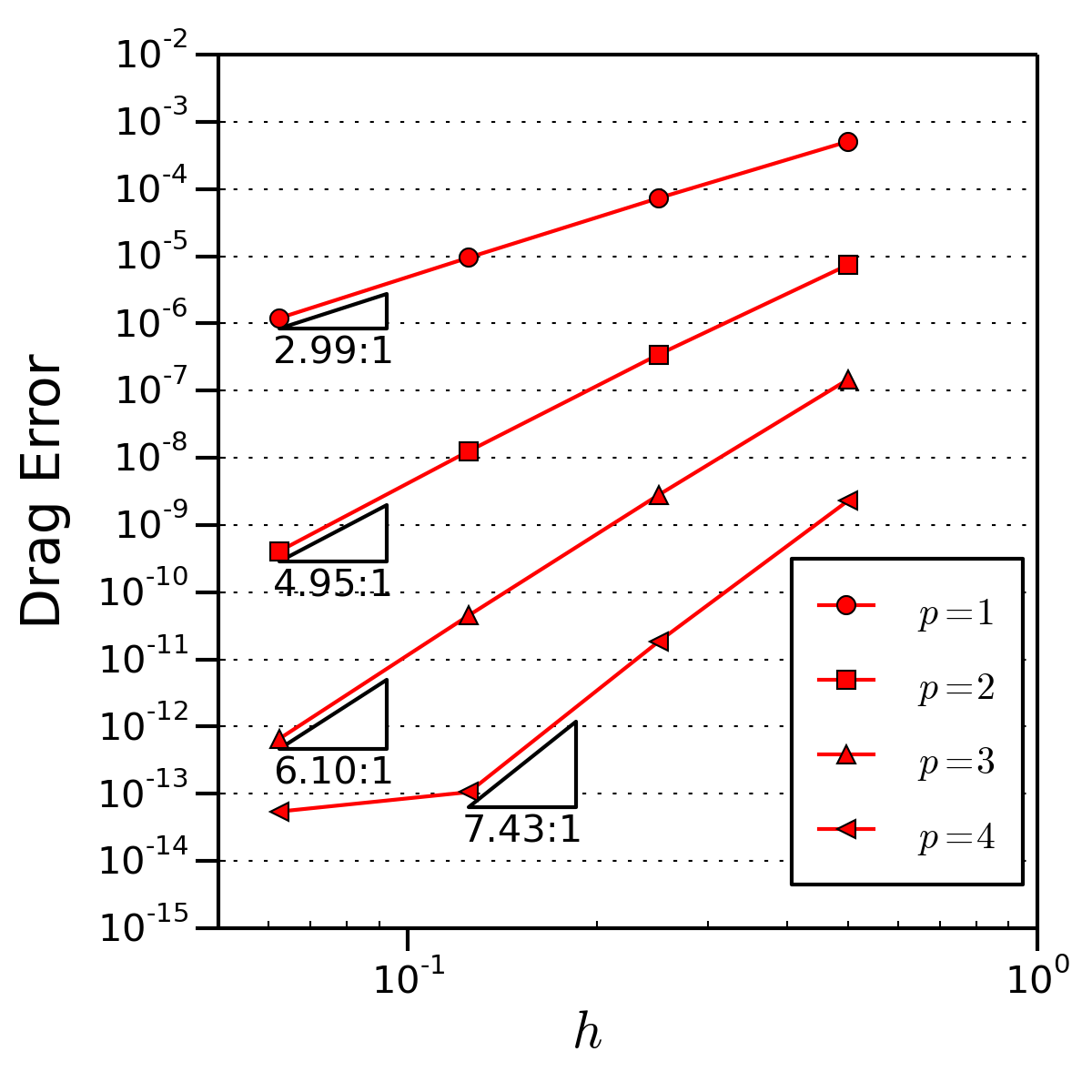}}
    \subfigure[DGD discretization\label{fig:dgd_func}]{%
    \includegraphics[width=0.49\textwidth]{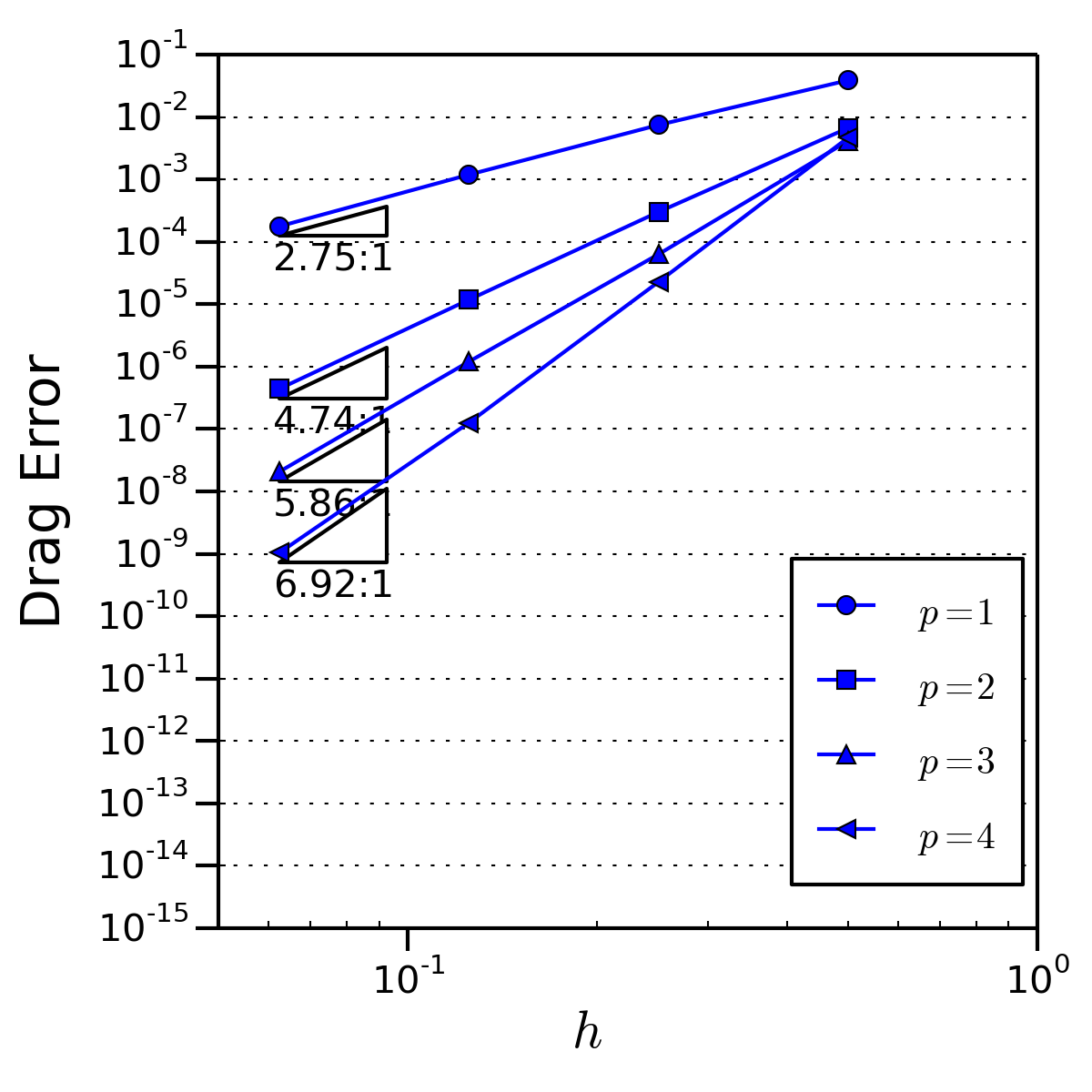}}
    \caption{Drag error of vortex exerted on inner radius}
    \label{fig:l2_results}
\end{figure}

We conclude our investigation of spatial accuracy by showing that the DGD discretization produces superconvergent functional estimates.  This study is motivated by the interest in integral outputs, such as lift and drag, in engineering analysis and design.  To this end, we estimate the drag force on the inner radius of $\Omega$ and compare this with the exact value of the drag.  Figure~\ref{fig:dsbp_func} shows that for degree $p = 1$ and $p=2$, the drag errors of the SBP scheme appear to approach a super-convergent rate of $2p+1$; for degrees $p=3$ and $p=4$, the convergence rates approach $2p$ and $2p-1$, respectively. Note that the drag error for $p=4$ stagnates around $10^{-13}$ due to the convergence tolerance used for Newton's method.  Similarly, Figure~\ref{fig:dgd_func} shows that the DGD drag errors are close to $2p+1$ super-convergent for degrees $p = 1$ and $p = 2$, while the convergence rates of the $p=3$ and $p=4$ schemes are closer to $2p$ and $2p-1$, respectively.

\subsection{Verification of entropy conservation/dissipation}

In this section, we apply the DGD discretization to the two-dimensional unsteady isentropic vortex problem to verify that total entropy is conserved, in the case of an entropy-conservative scheme, or dissipated, in the case of an entropy-stable scheme.  On an unbounded domain, the analytical solution to the unsteady vortex is given by 
\begin{equation} \label{eq:exact_solution2}
\begin{split}
    \rho &= \left(1-\frac{\epsilon^2 (\gamma-1)M^2}{8\pi^2}\exp{(f(x,y,t))}  \right)^{\frac{1}{\gamma-1}}, \\
    u &= \rho k \left(1-\frac{\epsilon}{2\pi} k(y - y_0) \exp{\left(\frac{f(x,y,t)}{2}\right)}\right), \\
    v &= \rho k^2 \frac{\epsilon}{2\pi}(x - x_0)\exp{\left(\frac{f(x,y,t)}{2}\right)}, \\
    e &= \frac{k^2p}{\gamma-1} + \frac{1}{2\rho}(u^2 +v^2),
\end{split}
\end{equation}
where $f(x,y,t) = 1 - ( (x-x_0-t)^2 + (y-y_0)^2)$, and the vortex is initially centered at $[x_0,y_0]^T = [0.5,0.5]^T$. The Mach number is $M = 0.5$, $\epsilon = 1.0$ is the vortex strength, and $k = 15$ is a scaling factor that controls the vortex speed.


We cannot simulate the vortex on an unbounded domain, so we consider the square domain $\Omega = \{(x,y)\;|\; x \in [0,1], y \in [0,1]\}$ with periodic boundary conditions applied on all four edges.  The analytical solution \eqref{eq:exact_solution2} with $t = 0$ defines the initial condition, but we do not make any claims regarding the convergence of the numerical solution to the exact solution, given the finite domain.  Despite possible solution errors, the total entropy should still behave as theoretically predicted, so this case remains useful for verification.

\begin{figure}[tbp]
\centering
    \subfigure[Example mesh\label{fig:unsteady_mesh}]{%
    \includegraphics[width=0.45\textwidth]{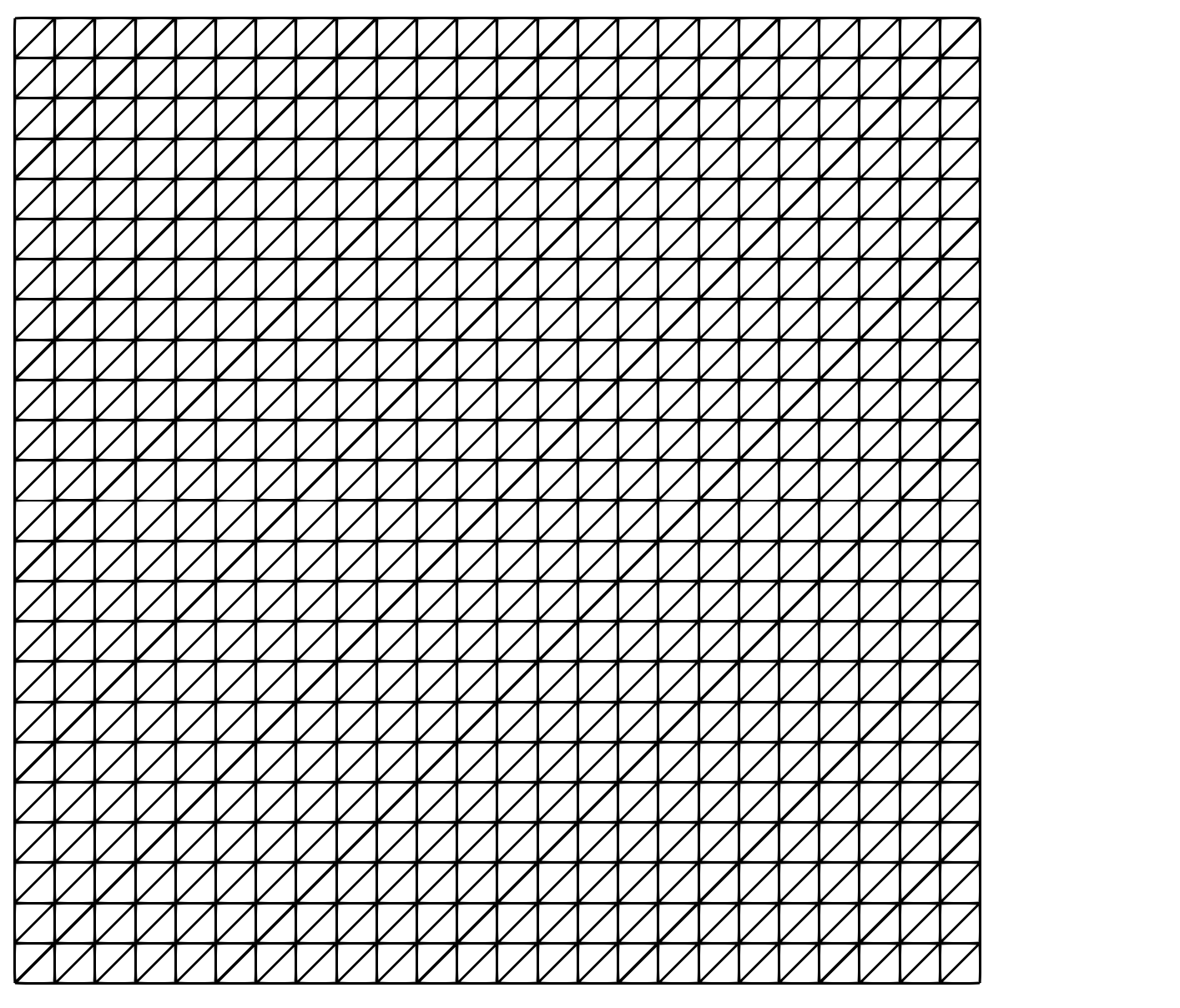}}
    \subfigure[Initial value of $\fnc{W}_1$ \label{fig:entropy_init}]{%
     \includegraphics[width=0.45\textwidth]{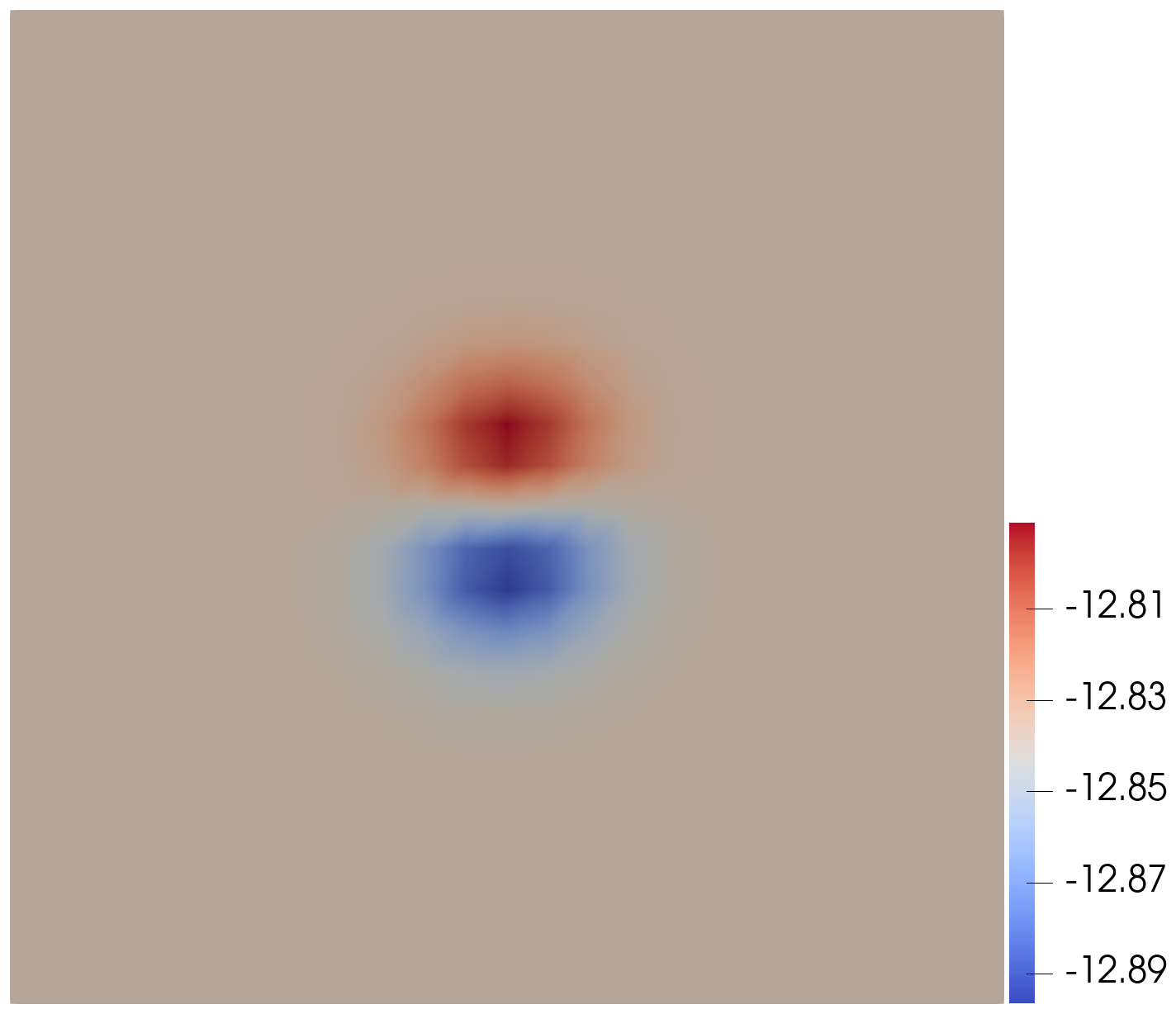}}
    \caption{Mesh and initial condition for the unsteady vortex study}
    \label{fig:unsteady_sample}
\end{figure}

All simulations are run for a period of time $T =1/k$ time units, during which the vortex travels in the positive $x$ direction and returns to its initial position.  The solution is advanced in time using the relaxation Runge-Kutta (RRK) \cite{Ranocha2020relaxation} variant of the implicit midpoint method described in Section~\ref{implicit_rrk}.  We use a constant CFL number of 10 based on the fastest acoustic wave speed, and the system \eqref{eq:DGD_midpoint} is solved to a tolerance of $1e-12$ using Newton's method.  Figure~\ref{fig:unsteady_sample} shows the mesh used for this study and the initial condition corresponding to the first entropy variable.

\begin{figure}[tbp]
\centering
    \subfigure[Entropy conservative scheme\label{fig:entropy_conservative_change}]{%
    \includegraphics[width=0.45\textwidth]{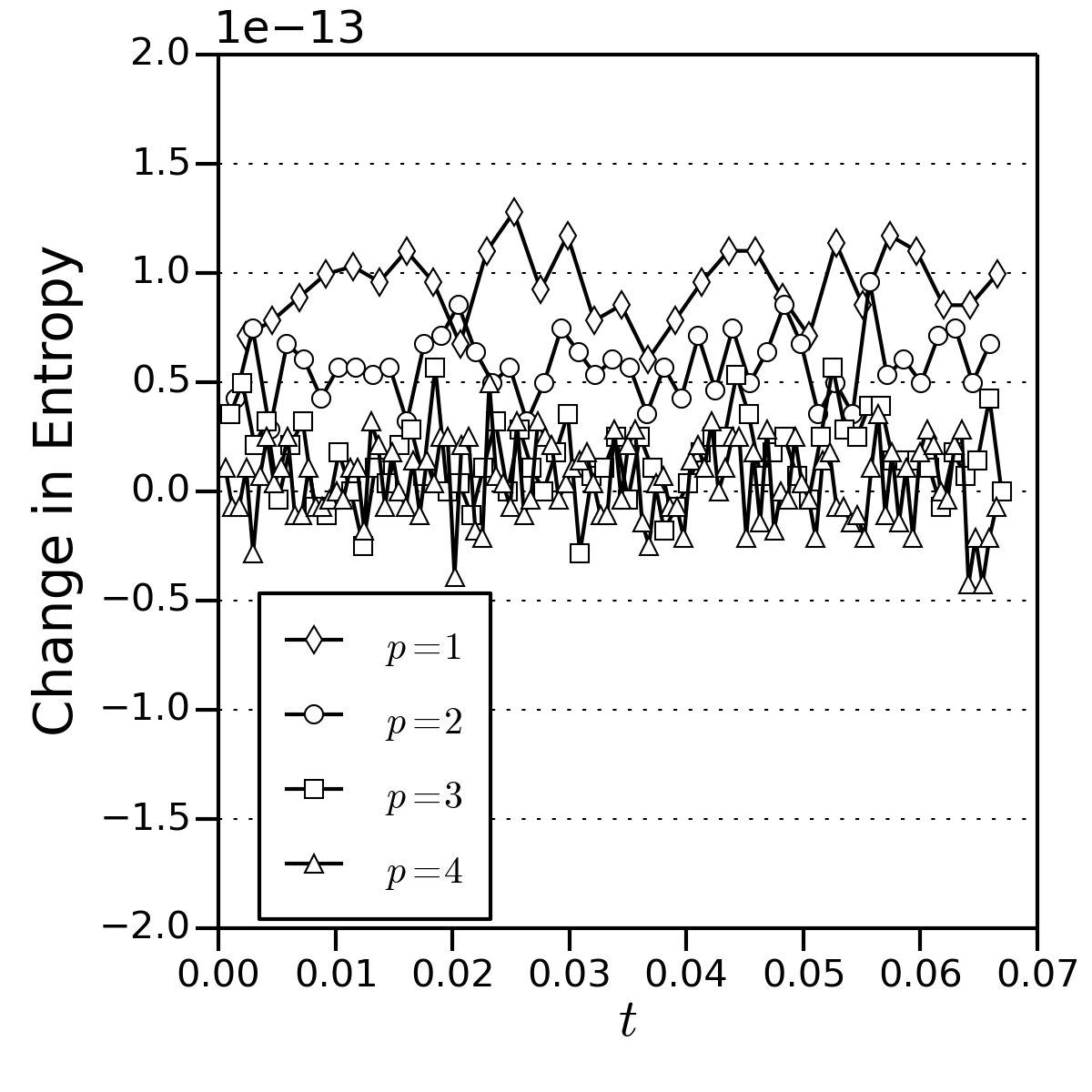}}
    \subfigure[Entropy stable scheme\label{fig:entropy_stable_change}]{%
     \includegraphics[width=0.45\textwidth]{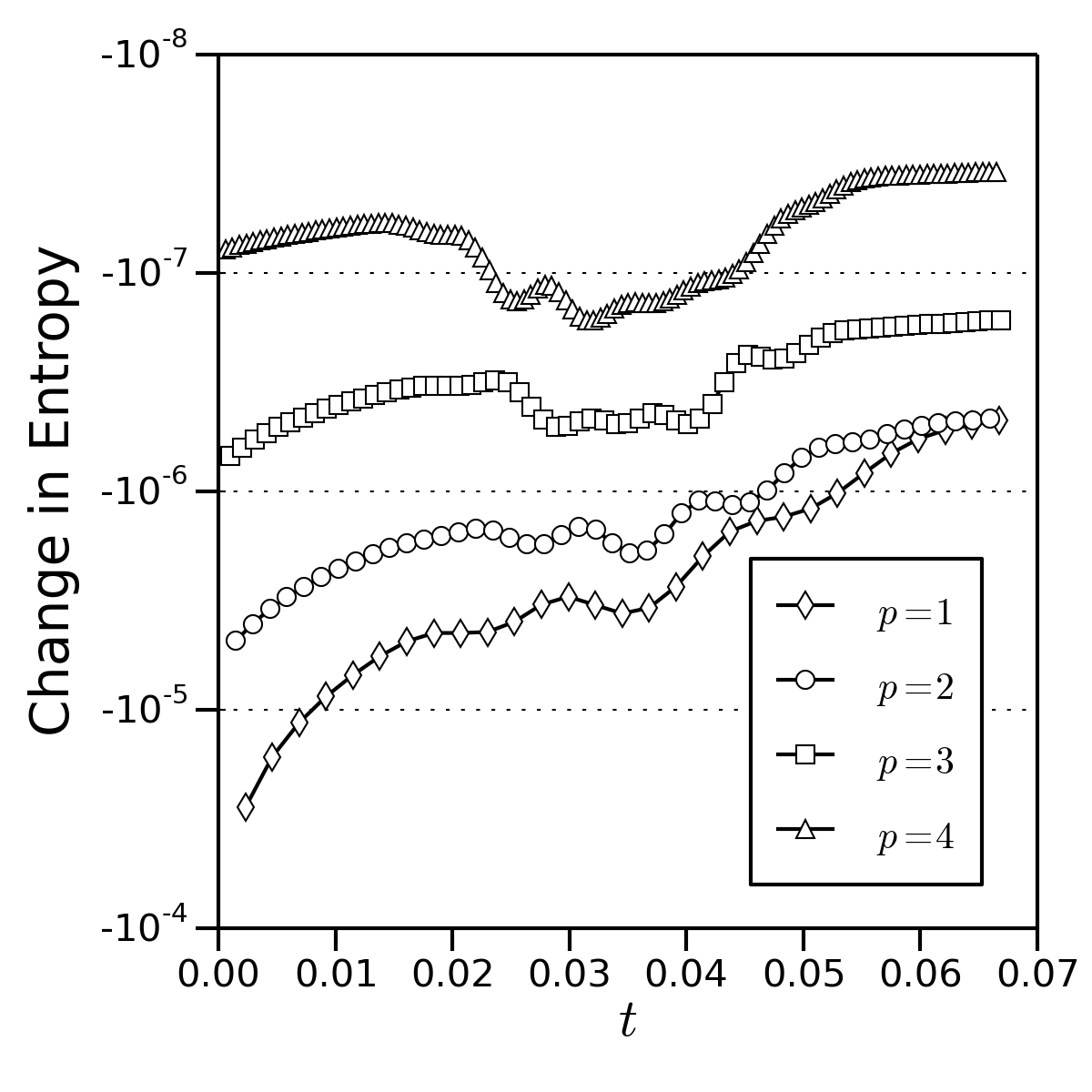}}
    \caption{Change in entropy for the unsteady vortex problem}
\end{figure}

Figure~\ref{fig:entropy_conservative_change} and Figure~\ref{fig:entropy_stable_change} show the time histories of the change in total entropy for the entropy-conservative and entropy-stable DGD schemes, respectively.  In both figures, the change in total entropy at time $t^{n+1}$ is given by
\begin{equation*}
 \Delta S^{n+1} = S(\wgd^{n+1}) - S(\wgd^{n}),
\end{equation*}
where $S(\wgd)$ is defined by~\eqref{eq:total_entropy}.  The entropy is conserved to order $10^{-12}$ for the entropy conservative scheme, which is consistent with the tolerances used in the Newton solver and secant method used for conserving the entropy among time steps.  Note that the piecewise-constant behavior in Figure~\ref{fig:entropy_conservative_change} is caused by the entropy changing in the last two digits only.  Finally, as predicted by the theory, the entropy change is always negative for the entropy-stable scheme.  Furthermore, less entropy is dissipated as the solution degree $p$ increases.

\subsection{Spectra of the DGD Jacobians}

This section investigates the spectra of the DGD Jacobians for the unsteady vortex problem.  Although the eigenvalues of the Jacobian are not relevant to the nonlinear entropy-stability of the DGD discretizations, the spectra reveal the potential efficiency of the schemes for conditionally stable time-marching methods and iterative solvers.  In addition, recent studies~\cite{Gassner2020stability,Ranocha2020preventing} indicate that some entropy-stable discretizations suffer from linear instabilities that raise doubts about their accuracy for long-time simulations.  The spectrum of the Jacobian will reveal if the entropy-stable DGD methods suffer from similar linear instabilities.

Linearizing the DGD semi-discretization~\eqref{eq:DGD_entstable} about the state $\wgd_0$, we obtain the following linear ODE for the perturbation $\wgd'$:
\begin{equation*}
    \frac{d \wgd'}{dt} = \tilde{\mat{A}} \wgd' + \bm{b},
\end{equation*}
where $\tilde{\mat{A}} \equiv -\tilde{\mat{M}}^{-1} \tilde{\mat{J}}$, and the global Jacobians are given by
\begin{equation*}
    \tilde{\mat{M}} = \sum_{k=1}^K \barPk^T \barMk \frac{\partial \uk}{\partial \wk} \barPk,
    \qquad\text{and}\qquad
    \tilde{\mat{J}} = \sum_{k=1}^K \barPk^T \frac{\partial \bm{r}_k}{\partial \wk} \barPk.
\end{equation*}
The element-level Jacobians $\partial \uk/\partial \wk$ and $\partial \bm{r}_k/\partial \wk$ are evaluated at the prolonged baseline state, $\barPk \wgd_{0}$.  The vector $\bm{b}$ is constant and is not relevant to the spectral analysis.

\begin{figure}[tbp]
\centering
    \subfigure[\label{fig:entropy_conservative_eigen}]{%
    \includegraphics[width=\textwidth]{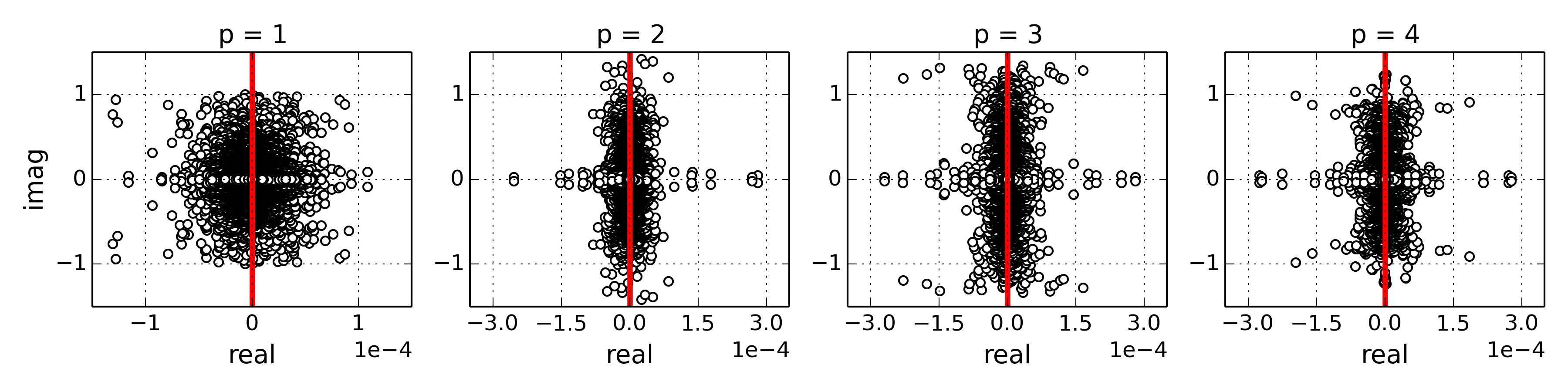}}
    \subfigure[\label{fig:entropy_stable_eigen}]{%
    \includegraphics[width=\textwidth]{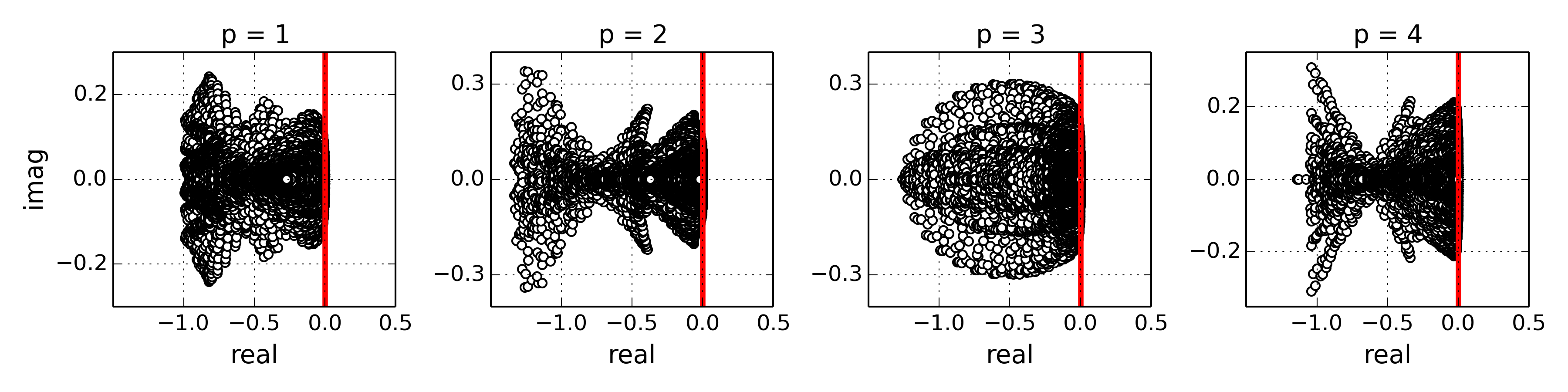}}
    \caption{Normalized eigenvalues from the entropy-conservative (upper row) and entropy stable (lower row) DGD discretizations}
\end{figure}

Figures~\ref{fig:entropy_conservative_eigen} and \ref{fig:entropy_stable_eigen} show the eigenvalues\footnote{We compute the eigenvalues of the generalized eigenvalue problem $\tilde{\mat{J}} \bm{v} = \lambda \tilde{\mat{M}} \bm{v}$} of $\tilde{\mat{A}}$ for the entropy-conservative and -stable schemes, respectively, run with 1152 triangular elements.  To normalize the results, the eigenvalues are scaled by the spectral radius of the corresponding $p=1$ operators.  The initial condition is adopted for the baseline state $\wgd_0$, which is used to evaluate the Jacobians $\tilde{\mat{M}}$ and $\tilde{\mat{J}}$; however, for the unsteady vortex, we do not expect the spectra to change significantly over time, because the exact solution is stationary under a suitable transformation.

It is noteworthy that the DGD spectral radius does not increase significantly as the polynomial degree increases.  This suggests that higher-order DGD schemes can take larger time steps relative to DG-type discretizations when using conditionally stable time-marching methods.  On the other hand, the appeal of using explicit time-marching schemes with the entropy-stable DGD schemes is diminished by the structure of their temporal terms; see Remark~\ref{remark:explicit}.

Finally, Table~\ref{tab:max_real} lists the maximum real part of the spectrum for each discretization on two mesh sizes, denoted by $h$ and $h/2$.  The table shows that all of the entropy-conservative discretizations are linearly unstable, with positive real parts on the order of $10^{-4}$.  The linear instability of the entropy-conservative schemes is not a significant concern, since some entropy dissipation is typically desirable for optimal convergence rates.  What is more concerning is that the entropy-stable schemes also appear to be linearly unstable, with the possible exception of a neutrally stable $p=1$ scheme.

\begin{remark} 
In the context of a (nonlinearly) entropy-stable scheme, linear instability does not necessarily imply that the solution will grow unbounded.  Nevertheless, linear instability may have implications for robustness --- since the growth can lead to a nonphysical state --- and for long-time accuracy.  For the DGD discretization of the unsteady isentropic vortex, the impact of the linear instability will take thousands of time steps to manifest itself, given the relatively small size of the positive real parts.
\end{remark}

\ignore{
Last but not least, the careful reader will notice that the real part of some eigenvalues in figure~\ref{fig:entropy_conservative_eigen} and \ref{fig:entropy_stable_eigen} is positive.  According to linear stability, this would imply that the discretizations are unstable.  However, this illustrates one of the limitations of linear theory and does not contradict nonlinear entropy stability theory.
}

\ignore{
\begin{table}[tbp]
    \centering
    \caption{Maximum real part over all eigenvalues in $\tilde{\mat{J}}$. \label{tab:J_max_real}}
    \begin{tabular}{lccccc}
    & & \multicolumn{4}{c}{\textbf{degree}} \\\cline{3-6}
    & \textbf{mesh size} & $p=1$ & $p=2$ & $p=3$ & $p=4$ \\\hline 
    \rule{0ex}{3ex}\textbf{entropy-conservative} & h &3.9154e-04  & 2.7077e-04 & 8.8772e-04 & 2.9911e-04 \\
    & h/2 &1.1685e-04 & 1.0139e-04 & 3.3225e-04 & 9.8273e-05\\\hline
    \rule{0ex}{3ex}\textbf{entropy-stable} & h &7.3513e-18  & 7.2077e-11 & 2.5203e-11 & 8.6476e-11 \\
    & h/2 & 6.7473e-17 & 1.9280e-11 & 2.7448e-11 & 1.2987e-11 \\\hline
    \end{tabular}
\end{table}
}

\begin{table}[tbp]
 \centering
 \caption{Maximum real part over all normalized eigenvalues in $\tilde{\mat{A}}$. \label{tab:max_real}}
 \begin{tabular}{lccccc}
 & & \multicolumn{4}{c}{\textbf{degree}} \\\cline{3-6}
 & \textbf{mesh size} & $p=1$ & $p=2$ & $p=3$ & $p=4$ \\\hline 
 \rule{0ex}{3ex}\textbf{entropy-conservative} & h &3.12e-04  & 2.76e-04 & 3.67e-04 & 3.33e-04 \\
 & h/2 &1.09e-04 & 2.81e-04 & 2.82e-04 & 2.77e-04\\\hline
 \rule{0ex}{3ex}\textbf{entropy-stable} & h &2.6909e-14  & 1.6337e-07 & 9.0194e-08 & 1.4249e-05 \\
 & h/2 &6.2607e-15 & 7.1420e-06 & 4.7396e-05 & 6.3008e-05 \\\hline
 \end{tabular}
\end{table}

\ignore{
\begin{table}[tbp]
 \centering
 \caption{Maximum real part over all eigenvalues in $\tilde{\mat{A}}$. \label{tab:max_real}}
 \begin{tabular}{lccccc}
 & & \multicolumn{4}{c}{\textbf{degree}} \\\cline{3-6}
 & \textbf{mesh size} & $p=1$ & $p=2$ & $p=3$ & $p=4$ \\\hline 
 \rule{0ex}{3ex}\textbf{entropy-conservative} & h &3.12e-04  & 2.76e-04 & 3.67e-04 & 3.33e-04 \\
 & h/2 &1.09e-04 & 2.81e-04 & 2.82e-04 & 2.77e-04\\\hline
 \rule{0ex}{3ex}\textbf{entropy-stable} & h &2.2774e-10  & 0.0014 & 7.6333e-04 &  0.1206 \\
 & h/2 &1.0715e-10 & 0.1222 & 0.8112 & 1.0784 \\\hline
 \end{tabular}
\end{table}
}

\subsection{Shock-tube problem}

We conclude the numerical experiments with a classical Riemann problem to investigate whether or not the DGD discretization correctly predicts shock speeds.  We are interested in this study because the conservative properties of the DGD discretization are not obvious: it has a nonlinear temporal term whose Jacobian is non-diagonal.

The Riemann problem is similar to the classical Sod's shock-tube problem.  The governing equations are the one dimensional compressible Euler equations and the space-time domain is given by $x \in [0,1], \; t \in [0, 0.3]$. The initial conditions are given by
\begin{equation*}
    \rho =     
    \begin{cases}
      5 \; & x < 0.5\\
      0.5 \; & x \geq 0.5
    \end{cases},
    \quad p =
    \begin{cases}
      1 \; & x < 0.5\\
      1/10 \; & x \geq 0.5
    \end{cases},
    \quad
    u = 0.
\end{equation*}
The initial density differs from the classical Sod problem.  This change was necessary because the present scheme does not include a shock-capturing method, and the oscillations at the shock produce  non-physical states when the classical initial condition is used.  

\begin{figure}[hbtp!]
\centering
    \includegraphics[width=0.8\textwidth]{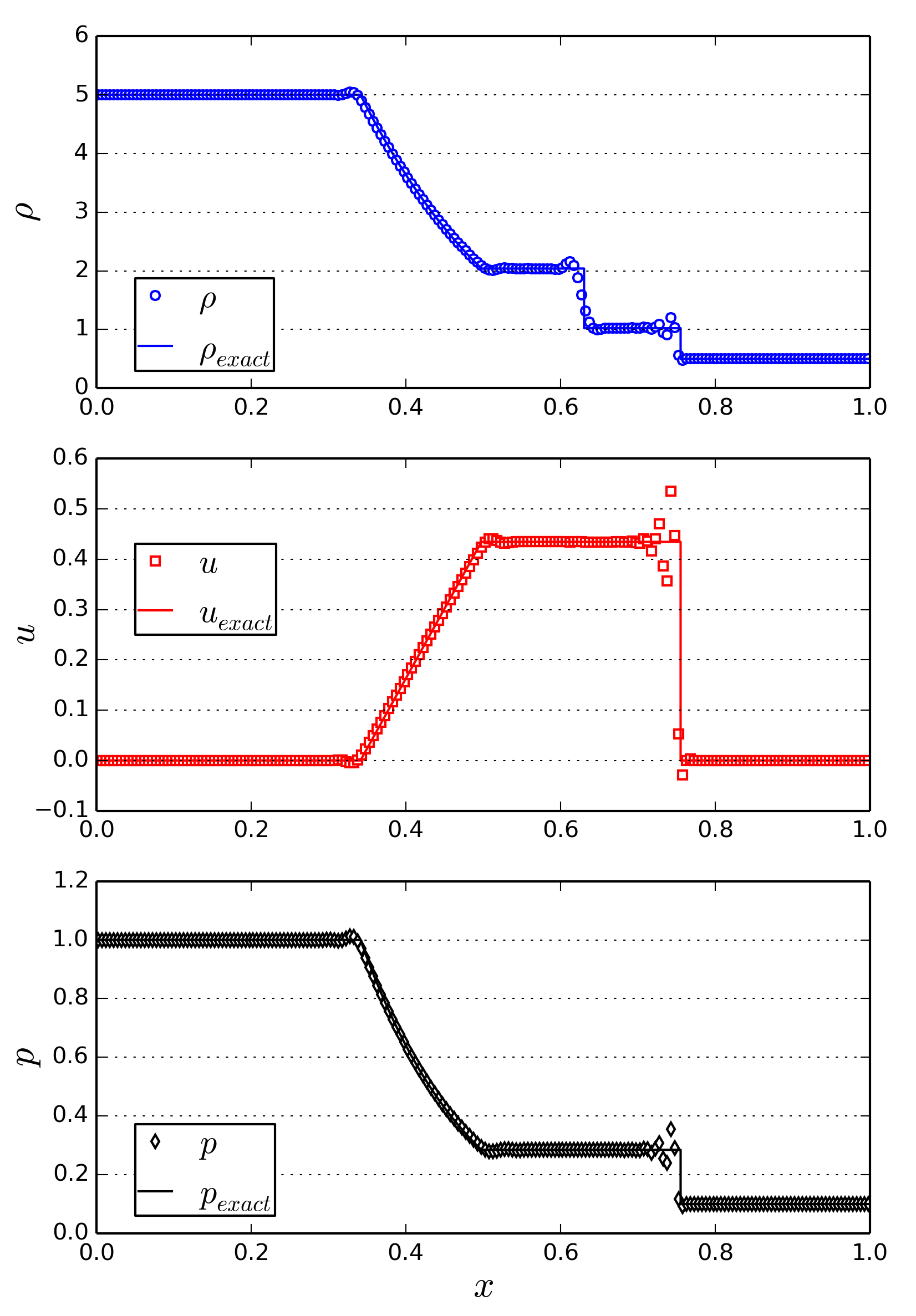}
    \caption{Exact and numerical solutions of Sod's shock problem.  The DGD solution ($p=2$) is plotted at the element centers.}
    \label{fig:sod_shock}
\end{figure}

Figure~\ref{fig:sod_shock} shows the density, velocity and pressure profiles at $t=0.3$ units based on the $p=2$ DGD discretization on a uniform mesh with 200 elements.  Oscillations are present at the shock and contact discontinuity, because, as explained above, the underlying scheme has no shock-capturing method.  However, the locations of the discontinuities agree with those of the exact solution. 

For a more quantitative assessment, Table~\ref{tab:sodshock_error} lists the $L^1$ errors in the first entropy variable, $\fnc{W}_1 = \partial \fnc{S}/\partial \rho$, for different mesh and degree configurations used to solve the shock-tube problem. We see that the $L^1$ error decreases with refinement, which suggests that the DGD scheme converges to the weak solution in the $L^1$ norm, despite the scheme’s nonlinear and non-diagonal temporal discretization.  Note that the rates of convergence for the even-order schemes are consistent with the results in~\cite{banks08a}.

\begin{table}[tbp]
\begin{center}
    \caption{$L^{1}$ errors for the Sod shock-tube problem, and estimated $L^{1}$ convergence rates. \label{tab:sodshock_error}}
    \begin{tabular}{lllllll}
    & & \multicolumn{2}{c}{$p=1$} & & \multicolumn{2}{c}{$p=2$}
    \\\cline{3-4}\cline{6-7}
    \rule{0ex}{3ex}$K$ & & 
    $L^1$ error & $L^1$ rate & & 
    $L^1$ error & $L^1$ rate \\\hline 
     \rule{0ex}{3ex}100 & &
     0.06436 & ---  & & 
     0.03333 & --- \\
     200 & &
     0.04116 & 0.6449 & & 
     0.02009 & 0.7303 \\
     400 & &
     0.02771 & 0.5708 & &
     0.01266 & 0.6662   \\
     800 & &
     0.02071 & 0.4201 & &
     0.008556 & 0.5653 \\\hline
     \rule{0ex}{6ex}& & \multicolumn{2}{c}{$p=3$} & & \multicolumn{2}{c}{$p=4$}
     \\\cline{3-4}\cline{6-7}
     \rule{0ex}{3ex}$K$ & & 
     $L^1$ error & $L^1$ rate & & 
     $L^1$ error & $L^1$ rate \\\hline 
     \rule{0ex}{3ex}100 & &
      0.03872 & ---  & & 
      0.02661 & --- \\
     200 & &
     0.02713 & 0.5132 & & 
     0.01599 & 0.7348 \\
     400 & &
     0.01968 & 0.4632 & &
     0.008951 & 0.8370  \\
     800 & &
     0.01380 & 0.5121  & &
     0.005142 &  0.7997 \\\hline
    \end{tabular}
\end{center}
\end{table}

\section{Summary and Conclusions}\label{conclusion}

We have presented entropy-conservative and entropy-stable discontinuous\linebreak Galerkin difference (DGD) discretizations of Euler equations on unstructured grids.  The entropy DGD discretization was constructed using three steps.
\begin{enumerate}
    \item Construct high-order prolongation/reconstruction operators to map DGD degrees of freedom from element centers to the nodes of appropriate degree diagonal-norm SBP operators.
    \item Use the prolonged DGD entropy variables in place of the SBP entropy variables in an entropy-conservative/stable SBP discretization.
    \item Apply the transposed prolongation operator to the SBP residual vector to map the SBP equations to the DGD equations.
\end{enumerate}

We verified the accuracy of DGD discretization by solving the steady vortex problem.  Furthermore, we demonstrated that the DGD and SBP discretizations have comparable efficiency when the $L^2$ error is measured versus number of degrees of freedom.  We also showed that the DGD discretizations produced superconvergent estimates of a drag functional.

The entropy-conservation and -stability properties were confirmed by solving the unsteady vortex problem. The entropy was conserved up to the accuracy of the iterative methods for the entropy-conservative scheme, and the mathematical entropy decreased monotonically for the entropy-stable schemes.

Our investigation of the spectra of the linearized entropy-stable DGD method revealed both positive and negative results.  First, the positive result: the spectral radius of the DGD discretization is relatively insensitive to order of accuracy, making it well suited to conditionally stable time-marching methods.  On the other hand, some eigenvalues had positive real parts for the high-order discretizations ($p \geq 2$), which indicates that the present entropy-stable DGD method is not immune to the linear instability observed in other entropy-stable SBP methods~\cite{Gassner2020stability,Ranocha2020preventing}.

Future work will include an investigation of conservation and shock capturing.  Our preliminary numerical experiments suggest that the DGD discretization is conservative and predicts shock speeds correctly; however, a formal proof that the scheme converges to the relevant weak solution has not been established.  In addition, shock capturing methods are needed to prevent non-physical states.

\section*{Acknowledgments}

G. Yan was supported by the National Science Foundation under Grant No. 1554253, and S. Kaur was supported by the National Science Foundation under Grant No. 1825991.  The authors gratefully acknowledge this support. We also thank RPI’s Scientific Computation Research Center for the use of
computer facilities.

\section*{Declaration}
\subsection*{Funding}
G. Yan and J. Hicken were supported by the National Science Foundation grant number No.1554253. S. Kaur was supported by the National Science Foundation grant number No.1825991
\subsection*{Conflicts of interest/Competing interests}
The authors declare that there have no conflict of interest.
\subsection*{Availability of data and material}
All data for the current study if available from the corresponding author upon request.
\subsection*{Code availability}
The software used to generate all data for the current study is available from the corresponding author upon request.

\appendix

\section{Proof of Lemma~\ref{lem:DGD_SBP_matrices}}\label{app:DGD_SBP_proof}

We begin with the right-hand side of the mass-matrix identity, and we recall that $(\mat{P}_k)_{qi} = \phi_i(\bm{x}_q)$:
\begin{equation*}
 \left[ \sum_{k=1}^{K} \mat{P}_k^T \Hk \mat{P}_k \right]_{ij} 
 = \sum_{k=1}^{K} \sum_{q=1}^{n_q} \phi_i(\bm{x}_q) (\Hk)_{qq} \phi_j(\bm{x}_q) = \sum_{k=1}^{K} \int_{\Omega_k} \phi_i  \phi_j \, d \Omega = \tM_{ij},
\end{equation*}
where we used the assumption that the diagonal entries in $\Hk$ and the nodes $X_k$ define a $2p$ exact quadrature for $\Omega_k$, and the fact that $\phi_i, \phi_j \in \mathbb{P}_p(\Omega_k)$ when restricted to element k.

For the $\tQx$ identity, we will use the SBP property \eqref{eq:SBP_accuracy}.  Specifically, for diagonal-norm SBP operators, \eqref{eq:SBP_accuracy} implies
\begin{equation*}
 \sum_{q=1}^{n_q} (\Qxk)_{rq} \phi_i(\bm{x}_q) = (\Hk)_{rr} \frac{\partial \phi_i}{\partial x}(\bm{x}_r),
    \qquad\forall\, \bm{x}_r \in X_k.
\end{equation*}
Note that the SBP operator differentiates $\phi_i$ exactly --- when restricted to element $k$ --- since this function is constructed from the degree $p$ polynomial basis $\{ \fnc{V}_j \}_{j=1}^{n_p}$.  Therefore, beginning with the right-hand side of the $\tQx$ identity, we find
\begin{align*}
\left[ \sum_{k=1}^K \mat{P}_k^T \Qxk \mat{P}_{k} \right]_{ij}
&= \sum_{k=1}^K \sum_{r=1}^{n_q} \sum_{q=1}^{n_q} \phi_i(\bm{x}_r) (\Qxk)_{rq} \phi_j(\bm{x}_q) \\
&= \sum_{k=1}^K \sum_{r=1}^{n_q}
\phi_i(\bm{x}_r) (\Hk)_{rr} \frac{\partial\phi_j}{\partial x}(\bm{x}_r) 
\\
&= \sum_{k=1}^K \int_{\Omega_k} \phi_i \frac{\partial \phi_j}{\partial x} \, d\Omega 
= (\tQx)_{ij}.
\end{align*}
In the penultimate step above, we applied the diagonal-norm quadrature to the degree $2p-1$ polynomial $\phi_i \frac{\partial \phi_j}{\partial x}$.

Finally, the $\tEx$ identity can be derived using the SBP operator property \eqref{eq:SBP_Ex_accuracy}:
\begin{align*}
\left[ \sum_{k=1}^K \mat{P}_k^T \Exk \mat{P}_{k} \right]_{ij}
&= \sum_{k=1}^K \sum_{r=1}^{n_q} \sum_{q=1}^{n_q} \phi_i(\bm{x}_r) (\Exk)_{rq} \phi_j(\bm{x}_q) \\
&= \sum_{k=1}^K \int_{\Omega_k} \phi_i \phi_j n_x \, d\Gamma = (\tEx)_{ij},
\end{align*}
where we again used the fact that the DGD basis functions are degree $p$ polynomials when restricted to an individual element.  This concludes the proof.

%


%
%

\bibliographystyle{spmpsci}      
\bibliography{jehicken}


\end{document}